\documentclass[12pt,leqno]{article}
\usepackage[dvipdfmx]{graphicx}
\usepackage{amsmath,amssymb}
\usepackage[pagewise]{lineno}
\usepackage{color}
\usepackage{here}
\usepackage{comment}

\newcommand{\red}{\textcolor[rgb]{1,0,0}}
\newcommand{\blue}{\textcolor[rgb]{0,0,1}}

\graphicspath{{./pic/}}

\makeatletter
 
 \@addtoreset{equation}{section}
\makeatother

\newcommand{\proof}{\noindent{\em Proof.} \;}
\newtheorem{thm}{Theorem}[section]
\newtheorem{lem}[thm]{Lemma}

\newtheorem{cor}[thm]{Corollary}
\newtheorem{prop}[thm]{Proposition}
\newtheorem{rmk}{Remark}[section]

\newcommand{\omu}{\overline{\mu}}

\newcommand{\oomega}{\overline{\omega}}
\newcommand{\ochi}{\overline{\chi}}

\newcommand{\bu}{\boldsymbol{u}}
\newcommand{\bL}{\boldsymbol{\Lambda}}
\newcommand{\bU}{\boldsymbol{U}}
\newcommand{\bV}{\boldsymbol{V}}
\newcommand{\brN}{\boldsymbol{N}}

\newcommand{\bR}{{\Bbb R}}

\renewcommand{\L}{{\mathcal L}}

\newcommand{\I}{{\mathcal I}}
\newcommand{\J}{{\mathcal J}}
\newcommand{\K}{{\mathcal K}}

\def\ep{{\epsilon}}

\def\la{{\langle}}
\def\ra{{\rangle}}

\def\eproof{\hfill{\vrule height5pt width3pt depth0pt}}

\begin{document}

\title{
Existence of  spiky stationary solutions to   \\
a mass-conserved reaction-diffusion model}
   

\author{
Yoshihisa Morita\\
\footnotesize Department of Applied Mathematics and Informatics,\\
\footnotesize Ryukoku University,\\
\footnotesize Seta Otsu 520-2194, Japan\\
\footnotesize({\tt morita@rins.ryukoku.ac.jp})\\
\\Yoshitaro Tanaka$^{*}$\\
\footnotesize Department of Complex and Intelligent Systems,\\
\footnotesize School of Systems Information Science,\\
\footnotesize Future University Hakodate,\\
\footnotesize 116-2 Kamedanakano-cho, Hakodate, Hokkaido 041-8655, Japan\\
\footnotesize({\tt yoshitaro.tanaka@gmail.com})
}

\date{}
\maketitle
\par\noindent


\begin{abstract}
We deal with a mass-conserved three-component reaction-diffusion system which 
is proposed by a model describing the dynamics of wavelike actin polymerization in the macropinocytosis and numerically exhibits dynamical patterns such as 
annihilation, crossover, and nucleation of pulses (Yochelis-Beta-Giv 2020).
In this article we first establish the condition for the diffusion driven instability in the system.
Then we rigorously prove the existence of spiky stationary solutions to the system
in a bounded interval with the Neumann condition.
By numerics these solutions play a crucial role in the nucleation of pulses. 
Reducing the stationary problem to a scalar second order nonlinear equation with a nonlocal term, we construct the desired solution by converting the equation to an integral equation.
\end{abstract}


\section{Introduction}
A certain class of reaction-diffusion systems are widely accepted 
as models describing pattern formations and dynamics in the fields of chemical reactions, morphogenesis, population biology, cell biology etc.
Those model equations have been attracting attention of mathematicians by providing rich mathematical problems concerned with variety of patterns and complex dynamical behaviors.

One of the central issues for studying such model equations is 
related to diffusion driven instability, so called Turing instability, which is a universal principle for the emergence of pattern formations. 
In many systems, once Turing instability takes place, that is, a spatially uniform steady state becomes unstable in the presence of diffusion, the system allows multi stable patterns.
However, there are model equations which exhibit a simple pattern though a Turing type instability occurs.
Such a model was proposed in \cite{OI} for understanding of
a mechanism of cell polarization by a mass-conserved reaction-diffusion 
system, where a mass constraint prevents complex patterns and
the stable steady state has a spatially simple profile.
Motivated by \cite{OI}, mathematical theory for stability and instability of nonuniform steady state and dynamical behaviors of solutions are developed in
\cite{MO}, \cite{M}, \cite{JM1}, \cite{JM2}, \cite{LS}, \cite{CMS} and \cite{LMS}.
It is also interesting that wave pinning of cell polarization is described
by a reaction-diffusion system with mass conservation \cite{MJK}.
Moreover, a model of asymmetric cell divisions is provided by mass conserved 
reaction-diffusion system in \cite{SLS} and a mathematical theory
is developed for the model in the work \cite{ML}.
In related to the dynamics for spiky pattern of mass-conserved reaction-diffusion system the readers may refer to \cite{KLE}. 
The readers may also refer to a nice review \cite{BHF} from a modeling and physical point of view.

In this article we are concerned with a new type of mass-conserved model
describing the dynamics of wavelike actin polymerization in  the process of macropinocytosis.
In this self-organization process circular dorsal ruffles (CDRs), which is dynamic ring-shaped undulations on the dorsal membrane, can be observed.
In order to understand the mechanism of wavelike dynamical process exhibited by CDRs, 
a reaction-diffusion system with a mass conservation was proposed (\cite{BDGY}).
Moreover, as a minimal model exhibiting qualitatively similar behaviors to those of 
the original system, the following three component reaction diffusion system was proposed (\cite{YBG}):

\begin{equation}
\label{eq:NSI}
\left\{
\begin{aligned}
&
N_t=\frac{N^2S}{1+I}-N+D_N\Delta N, \\
&
S_t=-\frac{N^2S}{1+I}+N+\Delta S,  \\
&
I_t=k_N N- k_I I+ D_I\Delta I,
\end{aligned}
\right. 
\end{equation}
where $N, S$ and $I$ stand for polymerized actin filaments
(F-actin), actin monomers (G-actin) and an actin polymerization inhibitor
respectively (see \cite{YBG} or \cite{BDGY} in details).
$\Delta$ stands for the Laplacian defined in a bounded domain $\Omega$ with the Neumann boundary condition.
We notice that in this model system a couple of nonlinear terms for F-actin and G-actin are expressed in  the form of switching, since the amount of protein in the cell is assumed to be
constant under no growth of the protein. 

As the main feature of the model equation \eqref{eq:NSI}, it is numerically shown that colliding pulse-like waves exhibit crossover and nucleation.  These behaviors also qualitatively correspond to those of CDRs in the cell membrane.   
The authors of \cite{YBG} address that this model can be distinguished from previous models by
robustness of the parameter range showing such a soliton-like behavior.

According to \cite{YBG}, we assume that the parameters satisfy
\[
0\le D_I \ll D_N <1,\qquad k_I<k_N.
\]
It is easy to see that the system allows a mass conservation as
\[
A:=\frac{1}{|\Omega|}\int_{\Omega} (N(x,0)+S(x,0))\ dx=\frac{1}{|\Omega|}\int_{\Omega} (N(x,t)+S(x,t))\ dx\quad(\forall t>0)
\] 
for appropriate initial data of \eqref{eq:NSI}.

It is challenging to reveal the dynamical structure of the model system rigorously and characterize soliton-like behaviors as reported in \cite{YBG}.
In this article, at the first step, we aim to examine the condition for the Truing-type instability rigorously and, in a bounded interval, prove that the system has a spiky equilibrium solution, which play a crucial role for the soliton-like behavior in the scenario demonstrated in \cite{YBG}.
In addition, we show that there exist other type of monotone solutions under a parameter condition.
In the next section we set the problem and state the main results of this paper. 

\section{Mathematical settings and main results}
We first introduce the new variables
and parameters as
\begin{align*}
&
u:=(k_N/k_I)N, \qquad v:=(k_N/k_I)S,\qquad  w:=I, \\
&
\kappa:=(k_N/k_I),\qquad \tau:=k_I,\qquad d=D_N, 
\qquad \ep:=D_I.
\end{align*}
Then the system \eqref{eq:NSI} in a bounded domain $\Omega\subset\bR^n$ reads
\begin{equation}
\label{uvw}
\begin{aligned}
&
u_t=d\Delta u+\frac{u^2}{\kappa^2(1+w)}v-u, \\
&
v_t=\Delta v-\frac{u^2}{\kappa^2(1+w)}v+u, \\
&
w_t=\ep \Delta w+\tau(-w+u)
\end{aligned}
\qquad \quad(x\in\Omega).
\end{equation}
We assume that the boundary $\partial\Omega$ is sufficiently smooth and set the Neumann boundary condition
\begin{align}
\label{Nbc}
\frac{\partial u}{\partial\boldsymbol{n}}=\frac{\partial v}{\partial\boldsymbol{n}}=
\frac{\partial w}{\partial\boldsymbol{n}}=0\qquad(x\in\partial\Omega),
\end{align} 
while the initial condition is given as
\begin{align*}
u(x,0)=u_0(x), \quad v(x,0)=v_0(x), \quad w(x,0)=w_0(x) \qquad (x\in \overline{\Omega}).
\end{align*}
We assume that $u_0(x), v_0(x)$ and $w_0(x)$ are taken so that there exists a unique global classical solution to \eqref{uvw} with \eqref{Nbc} (see \cite{R}). We don't discuss the condition for the global existence of \eqref{uvw} with \eqref{Nbc} which is out of the scope in the present article.

Since the mass conservation for the sum of $u$ and $v$ holds,  we thereby put
\begin{equation}
\label{M}
\begin{aligned}
M:=(k_N/k_I)A&=\langle u \rangle+\langle v\rangle\\
&:=\frac1{|\Omega|}\int_{\Omega}u(x,t)\ dx+\frac1{|\Omega|}\int_{\Omega}v(x,t)\ dx
\qquad(t\ge0).
\end{aligned}
\end{equation}

It is easy to see that \eqref{uvw} with \eqref{Nbc}
and \eqref{M} has a constant equilibrium solution $(u, v, w)=(0,M,0)$.
In addition, if 
\begin{align}
\label{Mc}
M>M_c:=\kappa^2+2\kappa,
\end{align}
then the system has other constant solutions 
\begin{equation}
\label{upm}
\begin{aligned}
(u,v,w)=&(u_{\pm}, v_{\pm}, u_{\pm}),\\
u_{\pm}:=&\frac1{2}(M-\kappa^2\pm\sqrt{(M-\kappa^2)^2-4\kappa^2}), \\
v_{\pm}:=&\frac1{2}(M+\kappa^2\mp\sqrt{(M-\kappa^2)^2-4\kappa^2}),
\end{aligned}
\end{equation}
where we used the condition $u+v=M$ for the constant solutions. 

We examine the linearized stability for each constant solution. In particular, $(u_+, v_+, u_+)$ allows the Turing-type instability in a parameter regime. 
Although it is not so easy to give an explicit condition on the parameters
for the Turing instability in the 3-component system, 
we can apply the result of \cite{ASY},
where two kinds of the instability for a general 3-component reaction diffusion systems are stated. 
The first one is called the steady instability (S-instability) whose definition is that one eigenvalue corresponding to
a nonzero wave number crosses $0$ from negative to positive along the real axis in the complex plane as a diagonal diffusion matrix suitably varies. 
The second one is called the wave instability (W-instability) which is the case
that non-real eigenvalue and its complex conjugate corresponding to a nonzero wavenumber cross the imaginary axis from the left to the right by changing diffusion matrix.
We have the following result:
\begin{prop}
\label{prop:T-inst}
For the the equilibrium $(u_+, v_+, u_+)$ of \eqref{uvw}, 
if $\tau \ge 1$, then the S-instability occurs 
for the diffusion matrix $D$
with $\max\{ d, \ep \} \ll 1$, 
while if $\tau < 1$,  both S-instability and W-instability occur for 
appropriate $D$ satisfying $\max\{ d, \ep \} \ll 1$.
\end{prop}

Next, we assume $\ep=0$ and $\Omega=(0,\ell)$ and
construct strictly monotone equilibrium solution to
\begin{equation}
\label{uv}
\begin{cases}
u_t=du_{xx}+f(u,v,w), \quad\\
v_t=v_{xx}-f(u,v,w), \quad \\
w_t=\tau(-w+u),
\end{cases}
\end{equation}
with the boundary condition 
\begin{align}\label{uvbc}
u_x=v_x=0\qquad(x=0, \ell),
\end{align}
and the same constraint \eqref{M}, where we put
\[
f(u,v,w):=\frac{u^2v}{\kappa^2(1+w)} -u. 
\]

As seen in \S\ref{proof_theorem}, we can reduce the stationary equations to a single equation
\begin{align}
\label{ugx}
du_{xx}+g(u;\mu)=0 \quad (0<x<\ell), \quad u_x=0\quad(x=0, \ell), 
\end{align}
with the constraint 
\begin{align}
\label{muM}
\mu:=M-(1-d)\la u \ra, \quad \la u \ra=\frac1{\ell}\int_0^\ell u(x)\ dx,
\end{align}
where we define
\begin{align*}
g(u;\mu)&:=u\left[\frac{u}{\kappa^2(u+1)}(\mu-du)-1\right], \qquad
G(u;\mu):=\int_0^u g(z;\mu)\ dz.
\\
\mu_c&:=\kappa^2+2\sqrt{d}\kappa.
\end{align*}
We note that $\mu_c$ is a critical value as 
the equation $g(u;\mu)=0$ has three distinct solutions for $\mu>\mu_c$.
Indeed, it
is easy to verify that if $\mu>\mu_c$, then, in addition to $u=0$,  $g(u;\mu)$ has two positive zeros
$\alpha(\mu)$ and $\beta(\mu)$, i.e.,
\begin{align*}
g(\alpha(\mu);\mu)=g(\beta(\mu);\mu)=0, \quad 0<\alpha(\mu)<\beta(\mu).
\end{align*}
Hence, in the parameter regime $\mu>\mu_c$ 
it is expected that the equation \eqref{ugx} exhibits a similar feature of 
the equation with a cubic nonlinearity. 
We have the following lemma:

\begin{lem}
\label{lem:1}
Assume $0<d<1$.
There exists a unique $\mu=\omu>\mu_c$ satisfying $G(\beta(\omu),\omu)=0$,
$G(\beta(\mu),\mu)<0~(\mu\in(\mu_c, \omu))$ and $G(\beta(\mu),\mu)>0~(\mu>\omu)$.
Then strictly monotone solutions to \eqref{ugx} are classified as
\begin{enumerate}
\item[\rm{(i)}]
For $\mu\in(\mu_c, \omu)$,  there exists $\ell_1(\mu)$ such that if $\ell>\ell_1(\mu)$, then 
\eqref{ugx} has a positive solution $u(\cdot; \mu, \ell)$ with $u_x(\cdot; \mu, \ell)>0~(0<x<\ell)$. Moreover, on any bounded interval in $[0,\infty)$, as $\ell\to\infty$,  $u(\cdot; \mu, \ell)$ uniformly converges to the homoclinic solution of
\begin{align}
\label{ugb}
du_{xx}+g(u;\mu)=0\quad(x\in\bR), \quad \lim_{|x|\to\infty}u(x)=\beta(\mu),
\quad u_x(0)=0.
\end{align}
\item[\rm{(ii)}]
For $\mu=\omu$, there exists $\ell_2(\omu)$ such that if $\ell>\ell_2(\omu)$, then 
\eqref{ugx} has a positive solution $u(\cdot; \omu, \ell)$ with $u_x(\cdot; \omu, \ell)>0~(0<x<\ell)$. Moreover, on any bounded interval in $[0,\infty)$, as $\ell\to\infty$,  $u(\cdot +\ell/2; \omu, \ell)$ uniformly converges to the heteroclinic solution to
the equation 
\begin{equation}\label{sol:hetero}
du_{xx}+g(u;\omu)=0~~~(x\in\bR), \quad
\lim_{x\to-\infty}u(x)=0, ~~ \lim_{x\to\infty}u(x)=\beta(\omu).	
\end{equation}
\item[\rm{(iii)}]
For $\mu>\omu$, there exists $\ell_3(\mu)$ such that if $\ell>\ell_3(\mu)$, then 
\eqref{ugx} has a positive solution $u(\cdot; \mu, \ell)$ with $u_x(\cdot; \mu, \ell)<0~(0<x<\ell)$. 
Moreover, on any bounded interval in $[0,\infty)$, as $\ell\to\infty$,  $u(\cdot; \mu, \ell)$ uniformly converges to the homoclinic solution to 
the equation which is obtained by replacing the asymptotic behavior in \eqref{ugb} by $\lim_{|x|\to\infty}u(x)=0$.
\end{enumerate}
\end{lem}
The readers might suspect if we could construct a monotone increasing solution 
in the case (iii) of Lemma \ref{lem:1}, which seems to be consistent with the cases (i) and (ii). As stated in the lemma, however, we prove the convergence of
the solution as $\ell\to\infty$. It turns out to be easier to handle the monotone decreasing solution in $[0, \ell]$ since it locally uniformly converges to the 
homoclinic solution asymptotic to $u=0$ in $[0,\infty)$.

We let $u(\cdot; \mu, \ell)$ be a solution to \eqref{ugx} in the one of 
(i), (ii) and (iii) in Lemma \ref{lem:1}. Then a simple application of the strong
maximum principle to \eqref{ugx} yields
\[
\mu-du(x; \mu, \ell)>0~\qquad(x\in[0,\ell]),
\]
indeed, if the left hand side takes a nonpositive value, then it leads us to a contradiction by the strong maximum principle. 
Hence, 
\begin{align*}
(u,v,w)=(u(\cdot; \mu, \ell), \mu-du(\cdot; \mu, \ell), u(\cdot; \mu, \ell))
\end{align*}
gives a positive equilibrium solution to \eqref{uv} with \eqref{uvbc} but no constraint \eqref{M}.
In fact,  the mass $M$ is determined by the solution $u(\cdot; \mu, \ell)$ as 
\begin{align}
\label{Mmuu}
M=\mu+(1-d)\la u(\cdot; \mu, \ell) \ra.
\end{align}
In the sequel, 
\begin{prop}
\label{prop:sol}
Let $u(\cdot; \mu, \ell)$ be a solution to \eqref{ugx} in Lemma \ref{lem:1}.
Then there exists a positive equilibrium solution to \eqref{uv} with \eqref{uvbc} and
$M$ given by \eqref{Mmuu}.
\end{prop}
\begin{rmk}
{\em 
This result is not sufficient in the sense that it is unclear how we find the specific $M$ for a solution $u(\cdot;\mu,\ell)$, 
in fact, in order to determine $M$, we need to compute $\la u(\cdot;\mu,\ell)\ra$ which depends on not only $\mu$ but also $\ell$. 
We should discuss whether there could exist a positive equilibrium solution to
\eqref{uv} with \eqref{uvbc} and \eqref{M} for given $M>0$.
In other words, for given $M$, we look for $\mu$ and $\ell$ solving
\eqref{Mmuu}.
}
\end{rmk}
Our main result is as follows:
\begin{thm}
\label{thm:1}
Assume $0<d<1$.
Then arbitrarily given ${M}>\omu$ 
there exists $\ell_{M}>0$ such that for each $\ell>\ell_{M}$
the system \eqref{uv} with \eqref{uvbc} and \eqref{M} has a positive equilibrium solution $(u^*(x;\ell), v^*(x;\ell), w^*(x;\ell))$ satisfying
$u_x^*(x;\ell)<0~(0<x<\ell)$ and 
\begin{align*}
v^*(\cdot;\ell)=M-(1-d)\la u^*(\cdot;\ell)\ra-du^*(\cdot;\ell), \qquad w^*(\cdot;\ell)=u^*(\cdot;\ell). 
\end{align*}
Moreover, as $\ell\to\infty$, $u^*(\cdot;\ell)$ locally uniformly converges to $u^h(\cdot)$
in $[0, \infty)$, where $u^h(\cdot)$ 
is the homoclinic solution satisfying 
\[
du_{xx}+g(u;M)=0~~(-\infty<x<\infty), \quad u_x(0)=0,\quad \lim_{|x|\to\infty}u(x)=0.
\]
\end{thm}

\begin{rmk}
{\em
To obtain the solution $u^*(\cdot;\ell)$, we show that
$\la u(\cdot; \mu, \ell)\ra$ takes smaller as $\ell$ takes larger
uniformly for $\mu$ in some interval. Consequently, 
we have $\la u^*(\cdot;\ell)\ra\to0$ as $\ell\to\infty$.
}
\end{rmk}
\begin{rmk}
{\em
As seen in the statement, even though we set $M$ arbitrarily large,
the corresponding solution exists for sufficiently large $\ell$. 
On the other hand, as for the solutions in (i) and (ii) of Lemma \ref{lem:1}, 
the range of $M$ is restricted by $\mu\in(\mu_c, \omu]$ and \eqref{Mmuu}. 
Although we do not have a clear result as Theorem \ref{thm:1},
we show the convergence of
$\la u(\cdot;\mu,\ell)\ra$ as $\ell\to\infty$ for each fixed $\mu
\in(\mu_c, \omu]$ in Appendix. 
}
\end{rmk}
\begin{rmk}
{\em We set $\ep=0$ in the original system \eqref{uv}. 
This certainly makes our analysis simpler.
Besides this technical reason, for $\ep=0$, very similar qualitative dynamics to the case $0<\ep  \ll  1$ is observed as described in \cite{YBG}.
Thus, as the first step, our approach seems to be quite reasonable.
Even though $\ep=0$, to our best knowledges, the system is a new model that has not been mathematically studied yet.
}
\end{rmk}
\begin{rmk}
{\em 
Numerically, the solution of Theorem \ref{thm:1} seems to be
unstable, though we didn't check any parameter values.
Mathematically rigorous proof for the stability is a future project.
In related to the stability
the readers can refer to \cite{MO}, \cite{M}, \cite{JM1}, \cite{JM2} and \cite{LMS}
where the stability of nonconstant solutions to some two-component conserved reaction-diffusion systems is examined while in \cite{PS} the instability of standing pulse solutions on the whole space of one-dimension is shown,
though those results cannot directly apply to our system.
}
\end{rmk}

\begin{figure}[bt]
\begin{center}
\begin{tabular}{cc}
	\includegraphics[keepaspectratio, scale=0.45]{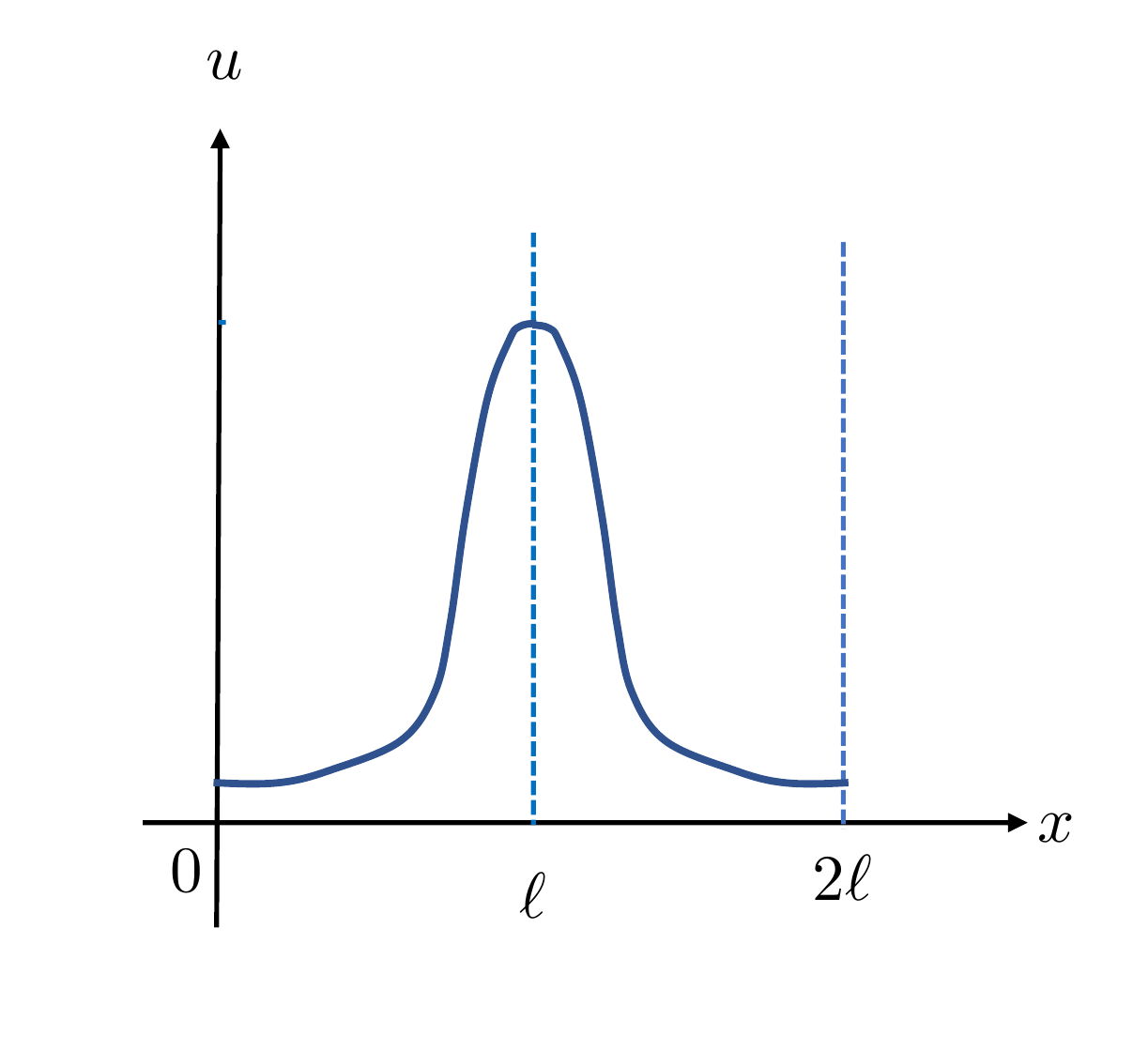}
	&\includegraphics[keepaspectratio, scale=0.45]{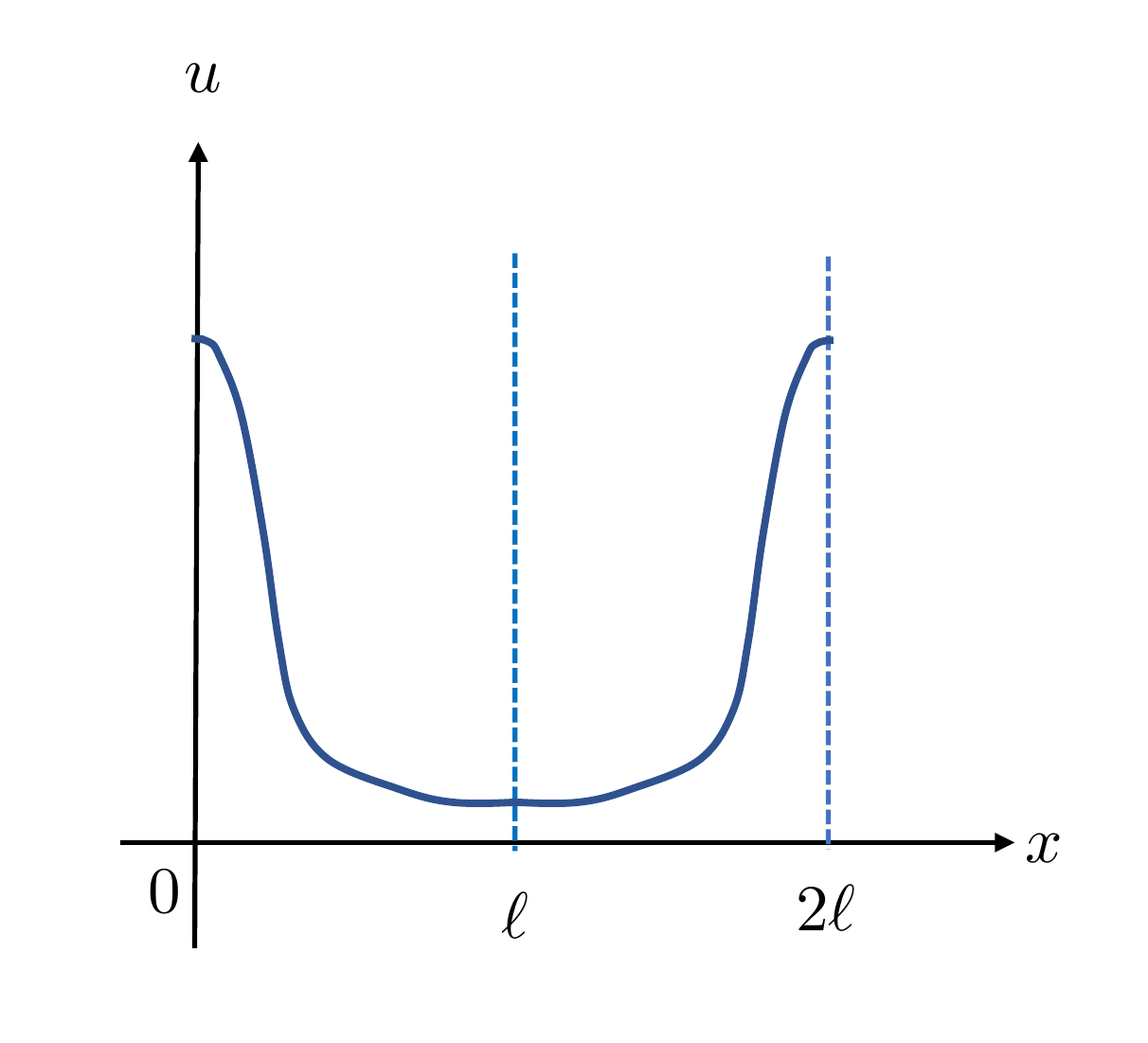} \\
	(a)&(b) \\
	\includegraphics[keepaspectratio, scale=0.4]{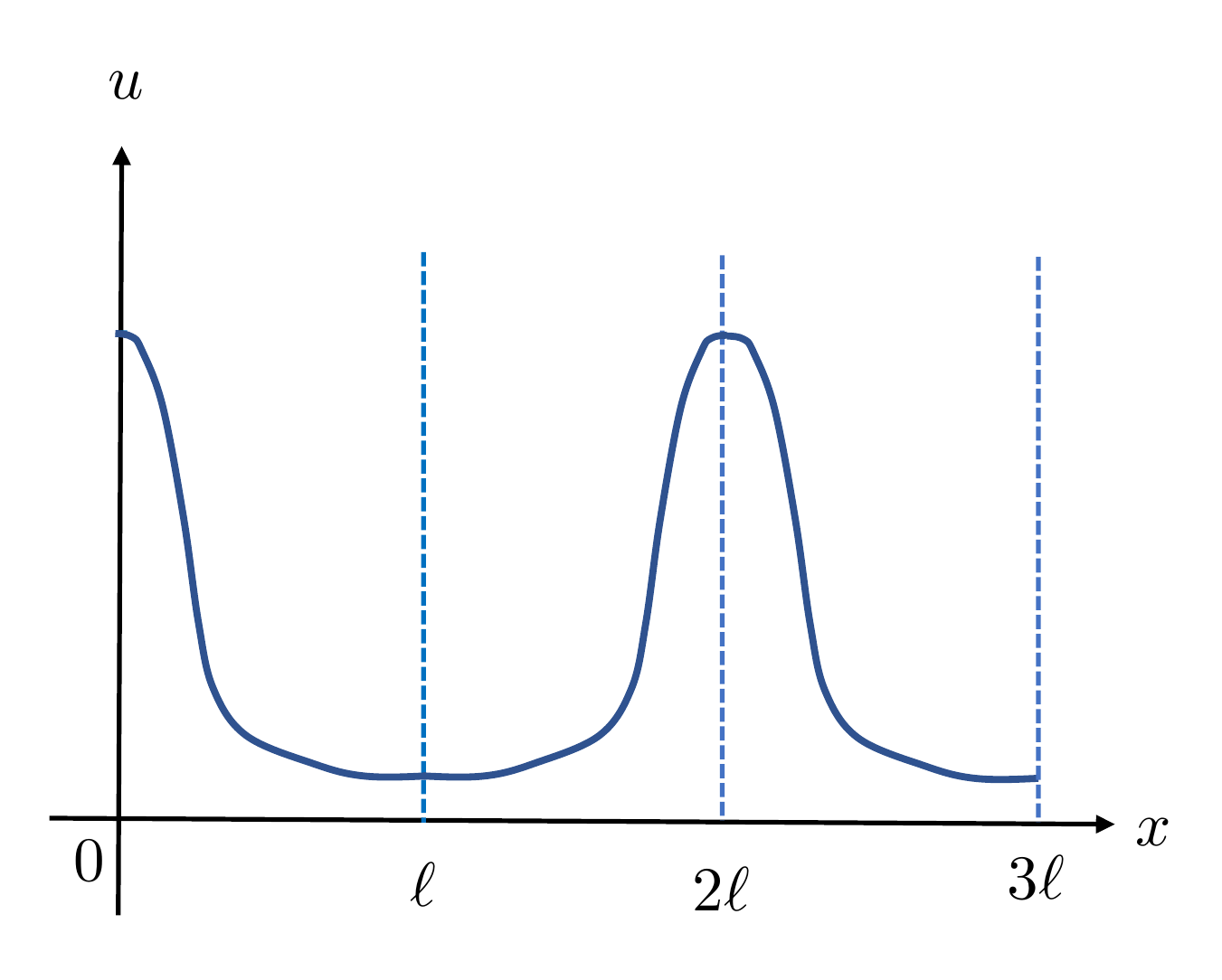}
	&\includegraphics[keepaspectratio, scale=0.4]{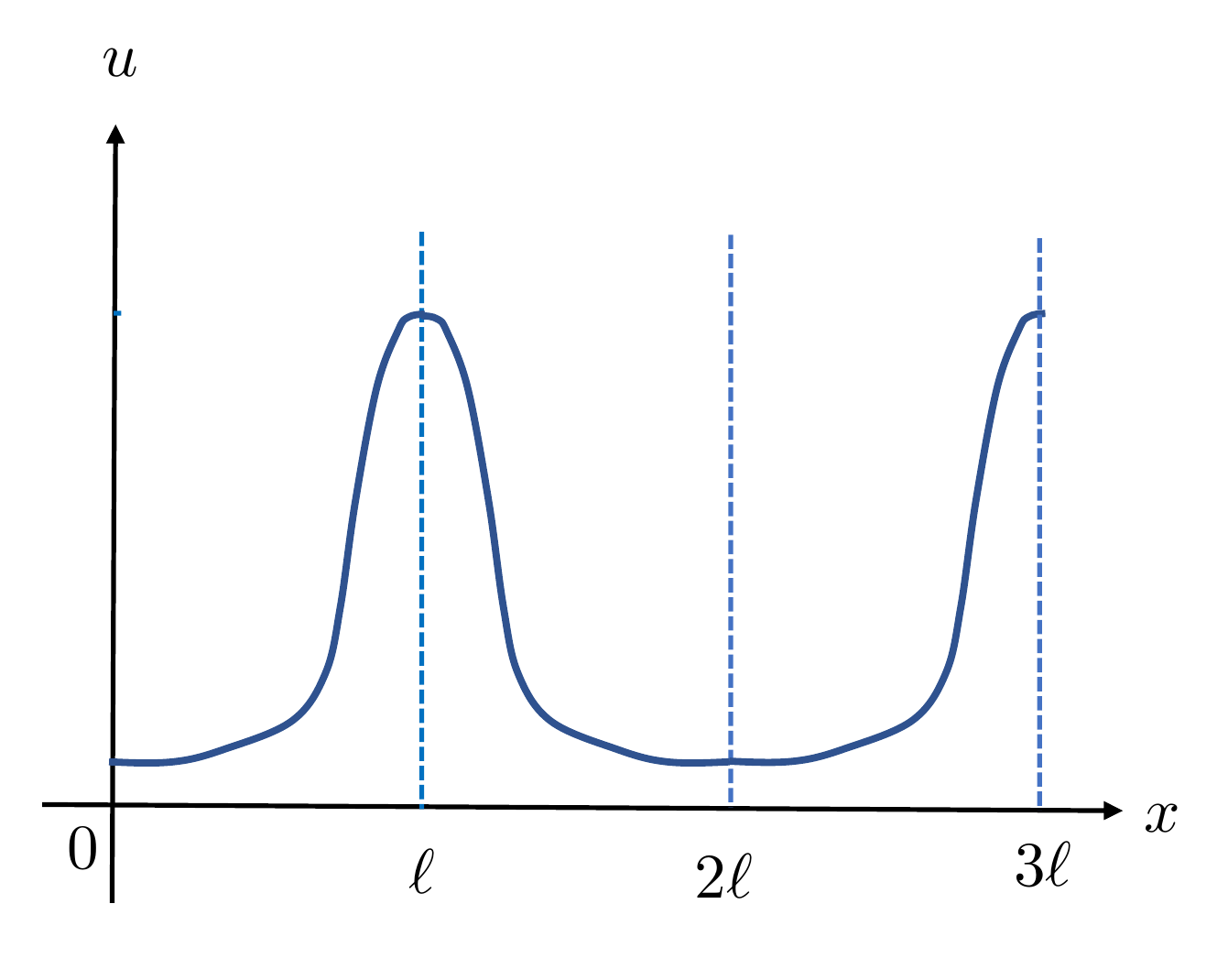} \\
	(c)&(d)
\end{tabular}
\end{center}
\caption{\small{
(a) Profile of the $u$-component of $\bL_1(\cdot;\ell)$. 
(b) Profile of the $u$ component of $\bV_1(\cdot;\ell)$.
(c) Profile of the $u$-component of $\bU_1(\cdot;\ell)$.
(d) Profile of the $u$-component of $\brN_1(\cdot;\ell)$.}}
\label{fig1}
\end{figure}

Making use of the solution $\bu^*(x;\ell)=(u^*(x;\ell), v^*(x;\ell), w^*(x;\ell))$ obtained in Theorem \ref{thm:1}, we can construct multi-mode solutions in  extended intervals. Define
\begin{align*}
&
\bL_1(x;\ell):=
\begin{cases} \bu^*(\ell-x;\ell) & (0\leq x\le \ell), \\
                      \bu^*(x-\ell;\ell) & (\ell\le x\le 2\ell), \end{cases}\\
&
\bV_1(x;\ell):=
\begin{cases} \bu^*(x;\ell) & (0\leq x\le \ell), \\
                      \bu^*(2\ell-x;\ell) & (\ell\le x\le 2\ell), \end{cases} \\
&
\bU_1(x;\ell):=
\begin{cases} \bu^*(x;\ell) & (0\leq x\le \ell), \\
                     \bL_1(x-\ell;\ell) & (\ell\le x\le 3\ell), \end{cases}\\
&
\brN_1(x;\ell):=
\begin{cases} \bu^*(\ell-x;\ell) & (0\leq x\le \ell), \\
                     \bV_1(x-\ell;\ell) & (\ell\le x\le 3\ell) \end{cases}
\end{align*}
(see Figure \ref{fig1})
and inductively for $j\geq2$ as
\begin{align*}
&
\bL_j(x;\ell):=
\begin{cases} \bL_{j-1}(x;\ell) & (0\leq x\le 2(j-1)\ell), \\
                      \bL_1(x-2(j-1)\ell;\ell) & (2(j-1)\ell\le x\le 2j\ell), \end{cases}\\
&
\bV_j(x;\ell):=
\begin{cases} \bV_{j-1}(x;\ell) & (0\leq x\le 2(j-1)\ell), \\
                      \bV_1(x-2(j-1)\ell;\ell) & (2(j-1)\ell\le x\le 2j\ell), \end{cases}\\
&
\bU_j(x;\ell):=
\begin{cases} \bU_{j-1}(x;\ell) & (0\leq x\le (2j-1)\ell), \\
                      \bL_1(x-(2j-1)\ell;\ell) & ((2j-1)\ell\le x\le (2j+1)\ell), \end{cases}\\
&
\brN_j(x;\ell):=
\begin{cases} \brN_{j-1}(x;\ell) & (0\leq x\le (2j-1)\ell), \\
                      \bV_1(x-(2j-1)\ell;\ell) & ((2j-1)\ell\le x\le (2j+1)\ell). \end{cases}
\end{align*}

%
Noticing 
\[
\la \bu^*(\cdot;\ell)\ra =(\la u^*(\cdot;\ell)\ra, \la v^*(\cdot;\ell)\ra, \la w^*(\cdot;\ell)\ra)=
\la \bL_j(\cdot)\ra=\la \bV_j(\cdot)\ra=\la \bU_j(\cdot)\ra=\la \brN_j(\cdot)\ra,
\]
we immediately obtain the next result.
\begin{cor}
Assume the same assumptions in Theorem \ref{thm:1} and $n_M:=[\ell/\ell_M]\geq 2$,
where $[\cdot]$ stands for the Gauss symbol.
If $n_M=2k$, then the system \eqref{uv} with \eqref{uvbc} and \eqref{M}
has equilibrium solutions $\bu=(u,v,w)$ given by
\begin{align*}
\bu=\begin{cases} 
            \bL_j(x;\ell'), ~~ \bV_j(x;\ell')\quad & (j=1, 2,\ldots, k)\quad \ell'=\ell/(2k),\\
          \bU_j(x;\ell''), ~~ \brN_j(x;\ell'')\quad  & (j=1, 2,\ldots, k-1), 
          \quad \ell'':=\ell/(2k-1),
          \end{cases} 
            \end{align*}
            while if $n_M=2k+1$, then the same assertion holds by replacing $k-1$ and $\ell/(2k-1)$ respectively by $k$ and $\ell/(2k+1)$ in the above. 
\end{cor}

\begin{rmk}
{\em When $g$ of \eqref{ugx} is cubic, global bifurcation structures were
extensively investigated under a similar constraint to 
\eqref{muM}  (see \cite{KMTY}, \cite{MKTY}, \cite{MTTY} and references therein). Their works are based on the expression for the solution by the elliptic integrals. On the other hand, in our problem we are not able to use the nice properties in the elliptic integrals as done in those works.}

\end{rmk}
In the next section we prove Proposition \ref{prop:T-inst}, and in \S\ref{proof_theorem} prove Theorem \ref{thm:1} together with the assertion (iii) of Lemma \ref{lem:1}.
In \S\ref{proof_lemma} we complete the proof of Lemma \ref{lem:1}.
Proposition \ref{prop:sol} clearly follows from Lemma \ref{lem:1}. 
In \S6 we give the proof of some lemmas in \S\ref{proof_theorem} and \S\ref{proof_lemma} which require lengthy computations.
We added an appendix where the convergence of $\la u(\cdot; \mu, \ell)\ra$ as $\ell\to\infty$ is shown for the solution $u(\cdot; \mu, \ell)$ in Proposition \ref{prop:sol}.

\section{Stability/instability of constant equilibria}
Nonnegative constant equilibria can be obtained by solving 
\[
\left(\frac{uv}{\kappa^2(1+u)}-1\right)u=0,\quad u+v=M.
\]
Hence, we have
\begin{align*}
& u=0 \qquad\mathrm{or}\nonumber \\ 
&
u(M-u)-\kappa^2(1+u)=-u^2+(M-\kappa^2)u-\kappa^2=0.
\end{align*}
This yields $(u, v, w)=(0, M, 0)$
and \eqref{upm} if \eqref{Mc} holds.
We note
\begin{align}
\label{ueq}
u_\pm^2=(M-\kappa^2)u_\pm-\kappa^2,
\end{align}
and
that when $M=M_c$, we have $u_{+}=u_{-}=\kappa$ and $v_{+}=v_{-}=\kappa^2+\kappa$.
Namely, a saddle-node bifurcation takes place at $M=M_c$.


\subsection{Linearized stability of constant equilibria}
We investigate the linearized stability around the constant equilibria.
Then the Jacob matrix of $(f(u,v,w), -f(u,v,w), \tau(u-w))^T$ reads
\begin{align*}
J(u,v,w):=\begin{pmatrix} f_u & f_v & f_w \\ -f_u & -f_v & -f_w \\
\tau & 0 & -\tau\end{pmatrix}
= \begin{pmatrix} \displaystyle\frac{2uv}{\kappa^2(1+w)} -1 & \displaystyle\frac{u^2}{\kappa^2(1+w)} & 
-\displaystyle\frac{u^2v}{\kappa^2(1+w)^2}
\\ -\displaystyle\frac{2uv}{\kappa^2(1+w)}+1
& -\displaystyle\frac{u^2}{\kappa^2(1+w)} & 
\displaystyle\frac{u^2v}{\kappa^2(1+w)^2} \\
\tau & 0 & -\tau \end{pmatrix}.
\end{align*}
At the constant equilibria
\begin{align*}
J(0,M,0)=& \begin{pmatrix}  -1 & 0  & 0 \\ 1 & 0 & 0 \\
\tau & 0 & -\tau \end{pmatrix}, \\
J(u_{\pm},v_{\pm}, u_{\pm})
=&
\begin{pmatrix}  1 & u_{\pm}/v_{\pm}  & -u_{\pm}/(1+u_{\pm}) 
\\  -1 & -u_{\pm}/v_{\pm}  & u_{\pm}/(1+u_{\pm})  \\
\tau & 0 & -\tau \end{pmatrix}.
\end{align*}
Indeed, apply
\[
\frac{u_{\pm}v_{\pm}}{\kappa^2(1+u_{\pm})}=1,
\]
to $J(u_{\pm},v_{\pm}, u_{\pm})$ and we obtain the above expression.

Let $\{\sigma_j\}_{j=1,2,\ldots}$ be the eigenvalues of $-\Delta$ with the Neumann boundary condition arranged in increasing order as
\[
\sigma_1=0<\sigma_2\leq \sigma_3\leq \cdots,
\]
and let $\{\varphi_j\}_{j=1,2,\ldots}$ be the corresponding eigenfunctions
to $\{\sigma_j\}_{j=1,2,\ldots}$ orthonormalized so that 
\[
(\varphi_i, \varphi_j)_{L^2}=\delta_{ij}
\]
holds.
Then the eigenvalue problem of the linearized operator
\[
\begin{pmatrix} d\Delta & 0 & 0 \\ 0 & \Delta & 0 \\ 0 & 0 & \ep\Delta
\end{pmatrix}
+ J(0,M,0)
\]
is reduced to the eigenvalue problem of the matrices
\[
A_j:=\begin{pmatrix}  -d\sigma_j -1 & 0  & 0 \\ 1 & -\sigma_j & 0 \\
\tau & 0 &  -\ep\sigma_j -\tau \end{pmatrix}, 
\quad j=1,2,\cdots .
\]
Similarly, 
\[
\begin{pmatrix} d\Delta & 0 & 0 \\ 0 & \Delta & 0 \\ 0 & 0 & \ep\Delta
\end{pmatrix}
+ J(u_{\pm}, v_{\pm}, u_{\pm})
\]
is decomposed into 
\begin{align}
\label{Bj}
B_j:=\begin{pmatrix}  -d\sigma_j +1 & u_{\pm}/v_{\pm}  & -u_{\pm}/(1+u_{\pm})  \\  -1 & -\sigma_j -u_{\pm}/v_{\pm} 
& u_{\pm}/(1+u_{\pm})  \\
\tau & 0 & -\ep\sigma_j -\tau \end{pmatrix}, 
\quad j=1,2,\cdots .
\end{align}
$A_1$ has three eigenvalues
\[
\lambda=0, ~ -1, ~ -\tau,
\]
while the remaining $A_j ~(j\geq2)$ have negative eigenvalues.
This implies $(0, M, 0)$  is linearly stable.

We check the eigenvalues of
\[
B_1=
\begin{pmatrix}  1 & u_{\pm}/v_{\pm}  & -u_{\pm}/(1+u_{\pm})  \\  -1 &-u_{\pm}/v_{\pm} 
& u_{\pm}/(1+u_{\pm})  \\
\tau & 0 & -\tau 
\end{pmatrix},
\]
and discuss the stability for spatially uniform perturbations, in other words, the stability of the
equilibria in the diffusion-free equations of 
\eqref{uvw}.

We compute $\mathrm{det}(B_1-\lambda I)$.
Put 
\[
\alpha_0:=u_{\pm}/v_{\pm},\quad \beta_0 :=u_{\pm}/(1+u_{\pm}).
\]
Then
\begin{equation*}
\begin{aligned}
\mathrm{det}(B_1-\lambda I)
 &
 =-\lambda\{\lambda^2+(\alpha_0-1+\tau)\lambda+
 \tau (\alpha_0-1+\beta_0)\}=0.
\end{aligned}
\end{equation*}
We check the signs of the coefficients.
\begin{align*}
\alpha_0-1+\beta_0&=\frac{u_{\pm}}{v_{\pm}}-1+\frac{u_{\pm}}{1+u_{\pm}}
=\frac{u_{\pm}}{v_{\pm}}-\frac{1}{1+u_{\pm}} \\
&=
\frac{u_{\pm}^2}{\kappa^2(1+u_{\pm})} -\frac{1}{1+u_{\pm}} 
=\frac{u_{\pm}^2-\kappa^2}{\kappa^2(1+u_{\pm})}.
\end{align*}
In view of \eqref{upm} and \eqref{ueq},
\begin{align*}
u_{\pm}^2&=(M-\kappa^2)u_{\pm}-\kappa^2=
\frac{2\kappa^2(M-\kappa^2)}{M-\kappa^2\mp\sqrt{(M-\kappa^2)^2-4\kappa^2}}-\kappa^2 \\
&=
\kappa^2\left\{\frac{2(M-\kappa^2)}{M-\kappa^2\mp\sqrt{(M-\kappa^2)^2-4\kappa^2}}-1\right\},
\end{align*}
thus,
\begin{align*}
u_{\pm}^2-\kappa^2&=
\kappa^2\left\{\frac{2(M-\kappa^2)}{M-\kappa^2\mp\sqrt{(M-\kappa^2)^2-4\kappa^2}}-2\right\}\\
&=
2\kappa^2\left\{\frac{\pm\sqrt{(M-\kappa^2)^2-4\kappa^2}}{M-\kappa^2\mp\sqrt{(M-\kappa^2)^2-4\kappa^2}}\right\}.
\end{align*}
This implies
\begin{align*}
\alpha_0-1+\beta_0 ~\begin{cases} >0 & (u=u_{+}), \\
<0 & (u=u_{-}).\end{cases}
\end{align*}
Hence, $(u_{-}, v_{-}, u_{-})$ is always unstable.
\par\medskip

We next go to the stability of $(u_{+}, v_{+}, u_{+})$.
Define $r(M):= u_{+}/v_{+} -1 + \tau$.

\begin{lem}
$(u_{+}, v_{+}, u_{+})$ is stable for spatially uniform perturbation if $M>M_c$
and
\begin{equation}
\label{k1}
 r(M_c)= \tau - \frac{\kappa }{\kappa + 1}=  k_I-\frac{k_N}{k_N+k_I}>0.
\end{equation}
On the other hand, if the inequality of \eqref{k1} is reverse, then
there is $M_*>M_c$ such that
$r(M_*)=0$ and $(u_{+}, v_{+}, u_{+})$ is stable (resp. unstable) for $M>M_*$
(resp. $\blue{M_c}<M<M_*$). 
\end{lem}
\proof
We show the monotonicity of $r(M)$ with respect to $M (>M_c)$.
Since
\begin{align*}
&
\frac{d}{dM}\left\{\frac{M(M+\kappa^2+\sqrt{(M-\kappa^2)^2-4\kappa^2})}{(M+1)}\right\} \\
&=
\frac{1}{(M+1)^2}(M+\kappa^2+\sqrt{(M-\kappa^2)^2-4\kappa^2}) \\
&
\qquad+\frac{M}{M+1}\left(1+\frac{M-\kappa^2}{\sqrt{(M-\kappa^2)^2-4\kappa^2}}\right)>0,
\end{align*}
so $r(M)$ is strictly increasing for $M>M_c$.
Next, we compute that
\begin{align*}
r(M)&
=\frac{u_{+}}{v_{+}}-1+\tau \\
&=\frac{M-\kappa^2+
\sqrt{(M-\kappa^2)^2-4\kappa^2}}{M+\kappa^2-\sqrt{(M-\kappa^2)^2-4\kappa^2}}-1+\tau \\
&
=
\frac{2M}{M+\kappa^2-\sqrt{(M-\kappa^2)^2-4\kappa^2}}-2+\tau \\
&=
\frac{M(M+\kappa^2+\sqrt{(M-\kappa^2)^2-4\kappa^2})}{2\kappa^2(M+1)}-2+\tau.
\end{align*}
When $M=M_c=\kappa^2+2\kappa$, 
\begin{align*}
r(M_c)&=\frac{(\kappa^2+2\kappa)(2\kappa^2+2\kappa)}{2\kappa^2(\kappa^2+2\kappa+1)}-2+\tau
=\frac{\kappa+2}{\kappa+1}-2+k_I\\
&=-\frac{\kappa}{\kappa+1}+k_I
=k_I-\frac{k_N}{k_N+k_I}.
\end{align*}
If \eqref{k1} is true, then there is a small $\delta>0$ such that
the equilibrium $(u_{+}, v_{+}, u_{+})$ is stable for $M\in(M_c, M_c+\delta)$ while
it is unstable for $M\in(M_c, M_c+\delta)$ if 
\begin{equation}
\label{k2}
k_I-\frac{k_N}{k_N+k_I}<0.
\end{equation}
 By putting $M_* = M_C + \delta$, the assertion of the lemma is true.
Specifically, 
$(u_{+}, v_{+}, u_{+})$ is stable for spatially uniform perturbation if $M>M_c$ and \eqref{k1} hold.
\eproof

\begin{rmk}
{\em 
When \eqref{k2} holds, a Hopf bifurcation takes place at $M=M_*$.
Some numerics shows that an unstable limit cycle exists for
$M>M_*$.
Note that when $k_I=0.3$ and $k_N=2$, we have
\[
k_I-\frac{k_N}{k_N+k_I}=0.3-\frac{2}{2.3}<0,
\]
for which \eqref{k2} is realized.
}
\end{rmk}
\subsection{Turing instability}
In this subsection we prove Proposition \ref{prop:T-inst}, that is, 
the Turing instability for $(u_{+}, v_{+}, u_{+})$. 
Define the diffusion matrix as
\begin{equation*}
	D:=
	\begin{pmatrix} d & 0 & 0 \\ 0 & 1 & 0 \\ 0 & 0 & \ep
	\end {pmatrix}.
\end{equation*}
Since the stability under the spatially uniform perturbations is assured by 
\[
\alpha_0-1+\beta_0>0, 
\qquad
r(M)>0,
\]
we assume that \eqref{k1} or $M>M_*$ with \eqref{k2}.

Recall \eqref{Bj} and consider 
\begin{align*}
&
\mathrm{det}(B_j-\lambda I) \\
 &\quad
 =-[\lambda^3+\{(d+1+\ep)\sigma_j-1+\alpha_0+\tau\}\lambda^2 \\
 &\qquad
 +\{(d+\ep d+\ep)\sigma_j^2+(-1+d\alpha_0 +\tau(d+1)+\ep(\alpha_0-1))\sigma_j+\tau(\alpha_0-1+\beta_0)\}\lambda \\
 &\qquad
 +\{\ep d\sigma_j^2+(\tau d+\ep(-1+d\alpha_0))\sigma_j
 +\tau(-1+d\alpha_0+\beta_0)\}\sigma_j] =0.
 \end{align*}
 We notice 
 \[
 -1+\beta_0=-1+\frac{u_{+}}{1+u_{+}}=-\frac{1}{1+u_{+}}<0.
\]

Since it is not easy to confirm an explicit parameter condition 
for the Turing instability by examining the above cubic equation, we apply the result of \cite{ASY}.
Following the theorem of \cite{ASY}, 
we set the complementary pair of subsystems of $J(u_+, v_+, u_+)$ as
\begin{equation*}
	J_2 := -\frac{u_+}{v_+}, \quad 
	J_{13} := 
	\begin{pmatrix}
	1 & -\dfrac{u_+}{1+u_+} \\
	\tau & -\tau
	\end{pmatrix}.
\end{equation*}
Since $J_2$ is negative, it is not concerned with the instability. 
On the other hand, in view of
\begin{equation*}
	\text{ tr} J_{13} = 1 -\tau, \quad \text{ det} J_{13} = \frac{-\tau }{ 1 + u_+ } < 0,
\end{equation*}
the assertions in (iii) and (iv) of Theorem 1.1 \cite{ASY} 
immediately follows. This completes the proof of Proposition \ref{prop:T-inst}.
\eproof


\begin{figure}[bt]
	\begin{center}
		\includegraphics[width=13cm, bb=0 0 814 231]{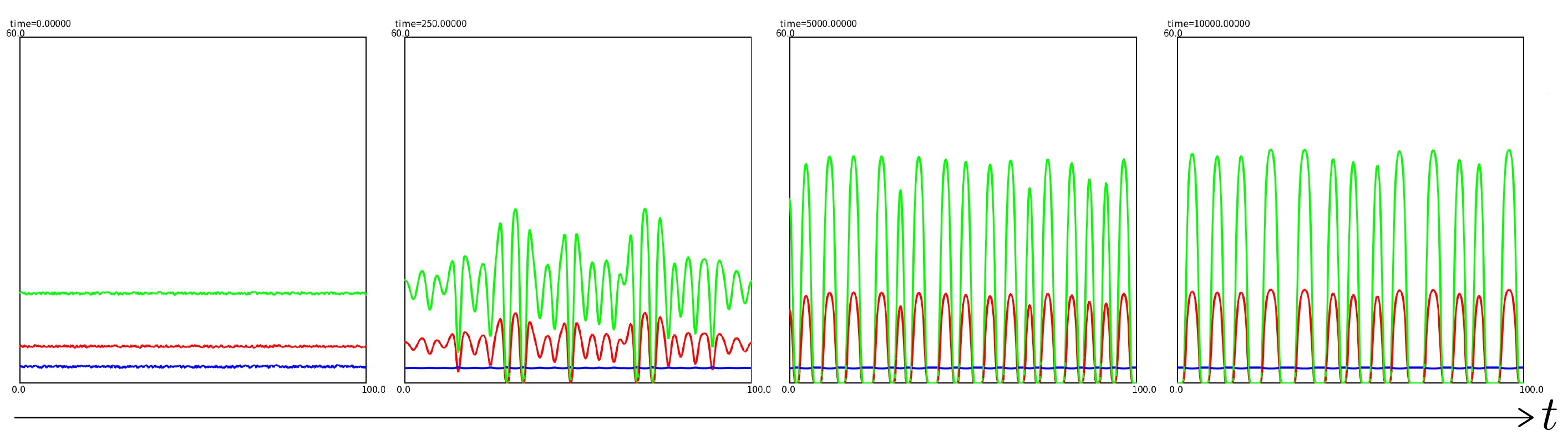}
	\end{center}
	\caption{\small{The results of a numerical simulation for the Turing instability.
			The parameters are $k_n=2.0, k_i=0.8, D_n=0.01,  D_i=0.001, A=8.8$. 
			The red, blue and green curves correspond to the solution of $N, S$ and $I$, respectively.
			The left, middle left, middle right and right pictures exhibit
			the profiles of solutions of \eqref{eq:NSI} in the interval $[0, 100]$ at $t=0, 250, 5000$ and $10000$, respectively. }}
	\label{fig2}
\end{figure}
In Figure \ref{fig2} numerically computed snap shots for $t=0, 250, 5000$ and $10000$ are displayed. Although we wanted to show the final state of the solution, 
even at $t=10000$, 
the wave pattern was not settled down yet.

\section{Proof of Theorem \ref{thm:1}}\label{proof_theorem}
We fist observe that the stationary problem of \eqref{uv} with 
\eqref{uvbc} and \eqref{M} turns to be
\begin{equation}
\label{equv}
\begin{aligned}
&
\begin{cases}
d\displaystyle{u_{xx}+u\left[\frac{uv}{\kappa^2(u+1)}-1\right]=0}, \\
\displaystyle{v_{xx}-u\left[\frac{uv}{\kappa^2(u+1)}-1\right]=0}
\end{cases}
\quad (0<x<\ell), \\
&
\qquad u_x=v_x=0\qquad (x=0, \ell), \\
&
\qquad
M=\la u\ra+\la v \ra,
\end{aligned}
\end{equation}
because of $u=w$.
Adding the two equations of \eqref{equv} yields
\[
du_{xx}+v_{xx}=0.
\]
By the Neumann condition this equation leads us to 
\[
du+v=d\la u\ra+\la v\ra.
\]
In the sequel, we obtain
\begin{equation}
\label{v}
v=-du+d\la u\ra +\la v\ra=-du+M-(1-d)\la u\ra.
\end{equation}
Plugging \eqref{v} into the first equation of \eqref{equv} 
yields \eqref{ugx} with \eqref{muM}. 
It is clear that by a solution $u^*(x)$ of \eqref{ugx},  $(u, v)=(u^*(x),\mu -du^*(x))$ gives a solution to  \eqref{equv} with no constraint $M=\la u\ra+\la v\ra$. 

We note that \eqref{ugx} has three constant solutions,
$u=0, \alpha(\mu)$ and $\beta(\mu)$, where
\begin{equation}
\label{ab}
\begin{aligned}
&
\alpha(\mu):=\frac1{2d}(\mu-\kappa^2-\sqrt{(\mu-\kappa^2)^2-4d\kappa^2}),\\
&
\beta(\mu):=\frac1{2d}(\mu-\kappa^2+\sqrt{(\mu-\kappa^2)^2-4d\kappa^2}),
\end{aligned}
\end{equation}
for $\mu>\mu_c=\kappa^2+2\sqrt{d}\kappa$.
By the definition of $\alpha(\mu)$ and $\beta(\mu)$, 
the equation of \eqref{ugx} is written as 
\begin{align}
\label{g:cubic}
du_{xx}+\frac{du(u-\alpha)(\beta-u)}{\kappa^2(u+1)}=0
\quad(0<x<\ell),
\end{align}
where we abbreviate $\alpha(\mu)$ and $\beta(\mu)$ as $\alpha$ and $\beta$ respectively. We will use this abbreviation as long as no confusion.

We emphasize that $\alpha(\mu)$ and $\beta(\mu)$ are not constant 
solutions to \eqref{ugx} with \eqref{muM} but  only to \eqref{ugx}. 
In other words, if there is $\mu$ satisfying
\[
\mu=M-(1-d)\alpha(\mu), \quad \text{or}\quad \mu=M-(1-d)\beta(\mu),
\]
then $\alpha(\mu)$ or $\beta(\mu)$  for such a $\mu$ gives a constant solution to \eqref{ugx} with \eqref{muM}, which must be one of $u_{\pm}$ defined in \eqref{upm}.

In view of  
\begin{align}
\label{Gexp}
\begin{aligned}
G(u;\mu)&=\frac{\mu}{\kappa^2}\left(\frac{u^2}{2}-u+\log(u+1)\right)\\
&\quad +
\frac{d}{\kappa^2}\left(-\frac{u^3}{3}+\frac{u^2}{2}-u+\log(u+1)\right)-\frac{u^2}{2},
\end{aligned}
\end{align}
and 
\begin{align}
\label{Gmu0}
G_\mu(u;\mu)=\frac{1}{\kappa^2}\left( \frac{u^2}{2}-u+\log(u+1) \right)>0\qquad(u>0),
\end{align}
we compute
\begin{align*}
\frac{d}{d\mu}G(\beta(\mu);\mu)&=g(\beta(\mu);\mu)\frac{d\beta}{d\mu}+G_{\mu}(\beta(\mu);\mu)>0.
\end{align*}
In addition, by $G(\beta(\mu_c),\mu_c)<0$ and the next lemma
there is a unique $\omu>\mu_c$ as stated in Lemma \ref{lem:1}.
%
\begin{lem}
\label{lem:mu}
Let
\begin{align}
\label{mu1}
\mu_1:=\kappa^2 +d+\frac{2}{3}\sqrt{ 3(4d\kappa^2+3d^2) }.
\end{align}
Then $G(\beta(\mu_1), \mu_1)>0$ holds.
\end{lem}
\proof
By a direct integration for $g(u;\mu)$ in the form of \eqref{g:cubic} we have
\begin{equation*}
	{\frac{\kappa^2}{d}}G(u;\mu) = -\frac{u^3}{3} + \frac{\alpha+\beta+1}{2}u^2+(\alpha+1)(\beta+1)(-u+\log(u+1)),
\end{equation*}
thus, 
\begin{equation*}
	{\frac{\kappa^2}{d}}G(\beta(\mu); \mu) = \frac{ \beta^3 }{6} -\frac{\alpha+1}{2}\beta^2 + (\alpha+1)  \Big\{ -\beta + (\beta+1)\log(\beta+1)  \Big\}
\end{equation*}
($\alpha=\alpha(\mu), \beta=\beta(\mu)$).
Since $(x+1) \log(x+1) -x >0$ for every $x\neq0$, we
have an estimate from below as
\begin{equation*}
	{\frac{\kappa^2}{d}}G(\beta(\mu); \mu) > \frac{ \beta(\mu)^3 }{6} -\frac{\alpha(\mu)+1}{2}\beta(\mu)^2.
\end{equation*}
It is easy to see that $  \beta(\mu) -3(\alpha(\mu)+1) \geq 0 $
implies $G(\beta(\mu);\mu)>0$.
Using \eqref{ab}, we solve the equation 
\begin{align*}
\beta(\mu)-3(\alpha(\mu)+1)=
\frac{2\sqrt{\mu^2 -(4d + 2\mu)\kappa^2 + \kappa^4}
+\kappa^2 - 3d - \mu  }{d}=0,
\end{align*} 
and obtain the solution $\mu=\mu_1$, which is as in \eqref{mu1}.
\eproof

\par\vspace{1cm}\noindent

Henceforth we assume $\mu>\omu$.
We put $u(0)=\xi>0$.
Multiplying \eqref{ugx} by $u_x$ and integrating it over $[0, x]$ 
leads us to
\begin{align*}
\frac1{2}u_x^2+G(u;\mu)=G(\xi; \mu),
\end{align*}
where we used $u_x(0)=0$.
Thus, a monotone decreasing solution is obtained by solving
\begin{align*}
\frac{du}{dx}=-\sqrt{2(G(\xi; \mu)-G(u;\mu))}.
\end{align*}

\begin{figure}[bt]
\begin{tabular}{cc}
\begin{minipage}[t]{0.45\hsize}
\centering
\includegraphics[keepaspectratio, scale=0.55]{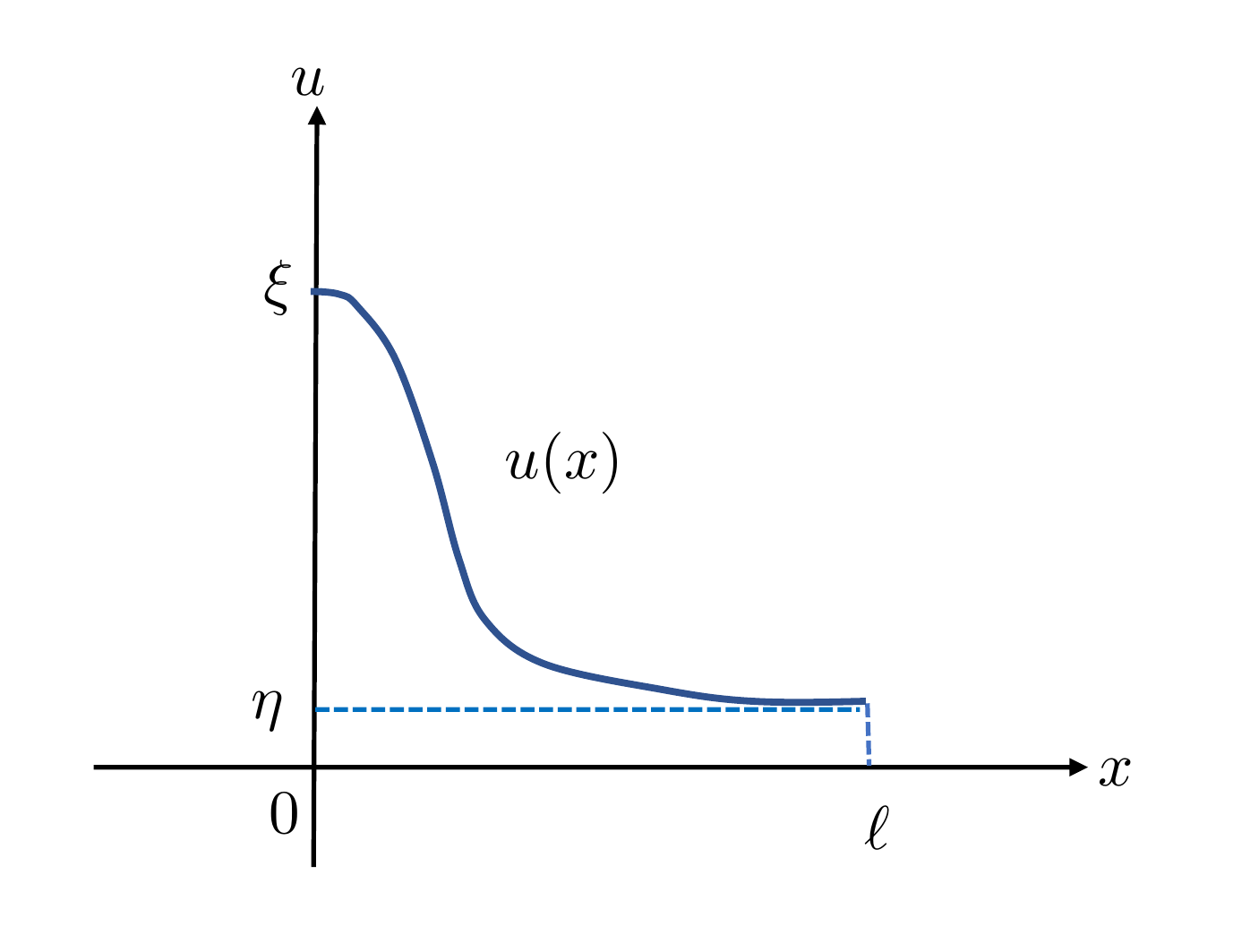}
\caption{\small{Profile of the solution $u(x)$}}
\label{fig3}
\end{minipage}
&
\begin{minipage}[t]{0.45\hsize}
\centering
\includegraphics[keepaspectratio, scale=0.55]{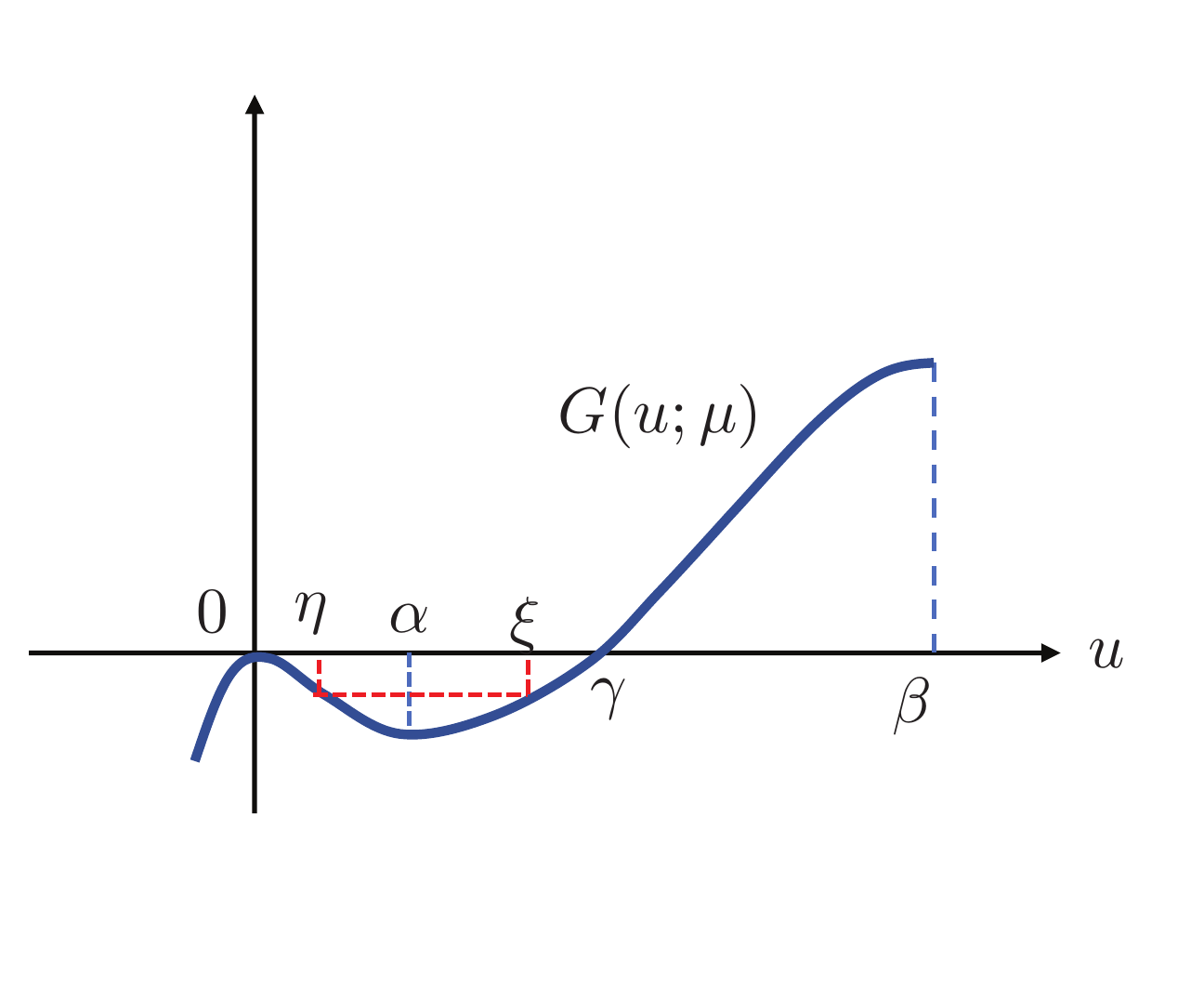}
\caption{\small{Graph of $G(u;\mu)$}}
\label{fig4}
\end{minipage}
\end{tabular}
\end{figure}

Integrating this yields
\begin{equation}
\label{xu}
\begin{aligned}
x&=\int_0^x 1\ dx \\
&
=-\int_{\xi}^{u}\frac{dz}{\sqrt{2(G(\xi; \mu)-G(z;\mu))}}
=\int_{u}^{\xi}\frac{dz}{\sqrt{2(G(\xi; \mu)-G(z;\mu))}}.
\end{aligned}
\end{equation}
A solution $u(x)$ satisfying 
\[
du_{xx}+g(u;\mu)=0\quad (0<x<\ell), \quad u(0)=\xi
\]
is obtained by the inverse function of \eqref{xu}. 
The existence of a monotone decreasing solution 
with $u_x(\ell)=0$, namely, a solution of \eqref{ugx} 
is achieved by $\xi$ and $\eta$ enjoying 
\begin{equation*}
\ell=\int_{\eta}^{\xi}
\frac{dz}{\sqrt{2(G(\xi; \mu)-G(z;\mu))}},
\end{equation*}
where $\eta$ is given by a function of $\xi$ satisfying $G(\xi; \mu) = G(\eta; \mu)$ 
(see Figures \ref{fig3} and \ref{fig4}).

In view of 
\[
G(\alpha, \mu)<0<G(\beta, \mu),\qquad
G_u(u;\mu)=g(u;\mu)>0\quad(\alpha<u<\beta),
\]
there is a unique zero $\gamma=\gamma(\mu)$ of
 $G(u;\mu)$ in $(\alpha, \beta)$, i.e.,
\begin{align}
\label{ga}
G(\gamma(\mu);\mu)=0, \quad \gamma(\mu)\in(\alpha, \beta),
\end{align}
by the intermediate value theorem
(see Fig.\ref{fig4}). 
Differentiating the equation \eqref{ga} with respect to $\mu$ yields
\[
\frac{\partial G(\gamma(\mu); \mu)}{\partial u}\frac{d\gamma(\mu)}{d\mu}
+\frac{\partial G(\gamma(\mu); \mu)}{\partial \mu}=0,
\]
which implies
\begin{align*}
\frac{d\gamma(\mu)}{d\mu}=-\frac{\partial G(\gamma(\mu); \mu)}{\partial \mu}\left/ \right.
\frac{\partial G(\gamma(\mu); \mu)}{\partial u}<0,
\end{align*}
where we used \eqref{Gmu0}
and
$\partial G(\gamma(\mu); \mu)/\partial u= g(\gamma(\mu); \mu)>0$.

\begin{lem}
\label{lem:eta}
There is a 
continuously differentiable function $\eta=\eta(\xi,\mu)\in(0,\alpha(\mu))$ 
in $(\xi, \mu)\in(\alpha(\mu),\gamma(\mu))
\times(\omu, \infty)$ satisfying
\begin{align}
\label{GG}
G(\eta(\xi,\mu) ;\mu)=G(\xi;\mu).
\end{align}
Moreover, $\eta(\xi,\mu)$ converges to $0$ as
$\xi\to\gamma(\mu)$ locally uniformly in $(\omu, \infty)$.
\end{lem}
\proof
Restriction of $G(u;\mu)$ on $[0, \alpha(\mu)]$
allows a continuous inverse function of it
because $y=G(u;\mu)$ is strictly monotone decreasing on $[0, \alpha(\mu)]$. We let $\Gamma(y; \mu)$ be 
such an inverse function. Then
we have an expression for $\eta$ as
\[
\eta(\xi,\mu)=\Gamma(G(\xi;\mu);\mu), \quad \xi\in[\alpha(\mu), \gamma(\mu)],
\]
with
$\eta(\alpha(\mu), \mu)=\alpha(\mu)$ and
$\eta(\gamma(\mu), \mu)=0$ by $G(\gamma(\mu);\mu)=0$. 
Since for 
$u\in(0,\alpha(\mu))$,
\[
G_u(u;\mu)=\frac{\partial }{\partial u}G(u;\mu)=g(u;\mu)<0
\]
holds, 
$\eta(\xi,\mu)$ is
continuously differentiable 
in $(\xi, \mu)\in(\alpha(\mu),\gamma(\mu))
\times(\omu, \infty)$.
In addition, by \eqref{Gexp}
\[
G(u;\mu)=-\frac{u^2}{2}+\frac{\mu}{3\kappa^2}(u^3+O(|u|^4))
+O(|u|^4)\qquad (|u|\ll1),
\]
in any bounded interval in $(\omu, \infty)$, 
we have 
\[
\Gamma(y;\mu)= -\sqrt{-2y}+O(|y|^3)
\]
holds in a small neighborhood of $y=0$ with $y<0$
for $\mu$ in the bounded interval.
Take $\mu_r>\omu$ arbitrarily and fix it.
Then for any $\ep>0$ there is $\delta_\ep>0$ such that
\begin{align}
\label{ee}
0<\sup_{\mu\in(\omu, \mu_r]}\eta(\xi,\mu)<\ep  \quad \text{holds}\quad \text{if} \quad \sup_{\mu\in(\omu, \mu_r]}|\xi-\gamma(\mu)|<\delta_\ep.
\end{align}
This concludes the proof.
\eproof
\par\medskip

By Lemma \ref{lem:eta}
the existence of the monotone decreasing solution to \eqref{ugx} is achieved by $\xi$ enjoying 
\begin{equation}
\label{Lxi}
\begin{aligned}
\ell=\rho(\xi, \mu):=&\int_{\eta(\xi, \mu)}^{\xi}
\frac{dz}{\sqrt{2(G(\xi; \mu)-G(z;\mu))}}\\
=
&\int_{\eta(\xi, \mu)}^{\xi}
\frac{dz}{\sqrt{2(G(\eta(\xi, \mu); \mu)-G(z;\mu))}}.
\end{aligned}
\end{equation}
\par\medskip
In order to prove that as $\ell\to\infty$, the desired solution converges to the homoclinic solution as stated in Theorem \ref{thm:1}, we have to make clear how the solution to \eqref{ugx} depends on $\mu$ and $\ell$.

\begin{rmk}
{\em
For each fixed $\mu>\omu$, one can easily prove $\lim_{\xi\to\gamma(\mu)}\rho(\xi,\mu)=\infty$ by modifying the argument below, so it implies the assertion (iii) of Lemma \ref{lem:1}.  We therefore skip the proof for  (iii) of Lemma \ref{lem:1}.
}
\end{rmk}

We take arbitrary $\mu_r>\omu$ and fix it below.
We show 
\begin{lem}
\label{lem:rho}
Arbitrarily given $N>0$, there is $\delta_N>0$, independent of $\mu\in(\omu,\mu_r]$, such that
for $\mu\in(\omu, \mu_r]$, 
\[
\rho(\xi,\mu)>N \quad \text{holds}\quad 
\text{if} \quad \xi\in(\gamma(\mu)-\delta_N, \gamma(\mu)).
\]
\end{lem}
We leave the proof of this lemma in the last part of this section.
\begin{lem}
\label{lem:dxi}
There exist $\xi_0(\mu)<\gamma(\mu)~(\mu\in(\omu, \mu_r])$
and $c>0$ such that 
\begin{align}
\label{rcM}
\frac{\partial\rho}{\partial \xi}(\xi, \mu)>c
\quad \text{in}\quad
\{(\xi, \mu): \xi_0(\mu)\leq \xi<\gamma(\mu), ~~ \omu<\mu\leq\mu_r\}
\end{align}
holds.
\end{lem}
Since the proof is lengthy, we will leave it the later section.

\par\medskip

By the above two lemmas we have the following result:
\begin{lem}
\label{lem:rell}
Let $\xi_0(\mu)$ be as in Lemma \ref{lem:dxi} and define
\[
\Sigma_r:=\{(\xi, \mu): \xi_0(\mu)\leq \xi<\gamma(\mu), ~~ \omu+\ep_r\leq\mu\leq\mu_r\}, 
\]
where $\ep_r$ is any small positive number satisfying $\omu+\ep_r<\mu_r$.
Consider the equation
\begin{align}
\label{rell}
\rho(\xi, \mu)=\ell,\qquad (\xi, \mu)\in \Sigma_r.
\end{align}
Then there is $\ell_0>0$ such that for each $\ell\geq\ell_0$, the equation \eqref{rell} has a unique solution $\xi=\xi(\mu, \ell)$ for each $(\mu, \ell)\in[\omu+\ep_r, \mu_r]\times[\ell_0,\infty)$. Moreover, $\xi(\mu,\ell)$ is 
$C^1$ in $\Sigma$ 
and satisfies 
\[
\frac{\partial\xi}{\partial\ell}(\mu, \ell)>0, \qquad
\lim_{\ell\to\infty}\xi(\mu,\ell)=\gamma(\mu).
\]
\end{lem}
\proof
Put $\ell_0:=\min_{(\xi, \mu)\in\Sigma}\rho(\xi, \mu)$.  Then the 
assertion of the lemma immediately follows from 
applying the implicit function theorem, Lemmas \ref{lem:rho} and \ref{lem:dxi}.
\eproof
\par\medskip

\par\medskip

By virtue of Lemma \ref{lem:rell}, for $(\mu, \ell)\in[\omu+\ep_r, \mu_r]\times[\ell_0,\infty)$,  we obtain a solution
$u=u(x; \mu, \ell)$ to \eqref{ugx} which satisfies $u(0; \mu,\ell)=\xi(\mu,\ell)$.

Moreover,
in \eqref{xu}, taking $\xi(\mu, \ell)\to\gamma(\mu)~(\ell\to\infty)$ yields
\begin{align*}
x=\int_{u}^{\gamma(\mu)}\frac{dz}{\sqrt{-G(z;\mu)}},
\end{align*}
and we can obtain an inverse function of this.
Then the reflection of the inverse function around $x=0$ leads to
a homoclinic solution $u^h(x; \mu)$ satisfying
\begin{align}
\label{uh}
du_{xx}+g(u;\mu)=0\quad(x\in\bR), \quad \lim_{|x|\to\infty}u(x)=0, 
\quad u_x(0)=0.
\end{align}
In the sequel, we obtain 
\begin{lem}
\label{lem:key}
Let $\ell_0$ be as in Lemma \ref{lem:rell}. 
Then for each $(\mu, \ell)\in[\omu+\ep_r, \mu_r]\times[\ell_0,\infty)$,
there is a positive and decreasing solution $u(\cdot;\mu, \ell)$ of \eqref{ugx}.
Moreover, in $[0, \infty)$, $u(\cdot;\mu, \ell)$ locally uniformly converges to the homoclinic solution $u^h(\cdot;\mu)$ 
satisfying \eqref{uh} as $\ell\to\infty$.
\end{lem}
\par\medskip\noindent
In addition, we have
\begin{lem}
\label{lem:avu}
For the solution $u(\cdot;\mu,\ell)$ obtained in Lemma \ref{lem:key}, 
$\int_0^\ell u(x;\mu, \ell)dx$ is continuous and uniformly bounded 
in $(\mu, \ell)\in [\omu+\ep_r, \mu_r]\times[\ell_0, \infty)$.
\end{lem}
\par\medskip\noindent
We leave the proof in the last part of this section.
The next lemma immediately follows from Lemma \ref{lem:avu}.
\begin{lem}
\label{lem:ur} 
Given $r>0$, there exists $\ell_r\geq \ell_0$ such that for $\ell\geq \ell_r$ 
\begin{align*}
\la u(\cdot;\mu,\ell) \ra=\frac{1}{\ell}\int_0^\ell u(x;\mu,\ell) dx\leq r
\qquad (\forall \mu\in[\omu+\ep_r, \mu_r])
\end{align*}
holds.
\end{lem}
Thus, we have
\begin{lem}
\label{lem:Lr}
Given $M>\omu$, 
take positive $r$ and $\ep_r$ so that $(M-\omu-\ep_r)/(1-d)>r$.
Then there exists $\ell_M>0$ such that for each 
$\ell\geq\ell_M$ the equation \eqref{ugx} with \eqref{muM} has
a monotone decreasing solution $u^*(x;\ell)>0$.
\end{lem}
\proof
It suffices to prove that there is $\mu\in[\omu+\ep_r, M]$ satisfying 
\eqref{muM}. 
Define 
\[
H_M(\mu):=\frac{M-\mu}{1-d},
\]
which is monotone decreasing in $\mu$.
Since $H_M(\omu+\ep_r)=(M-\omu-\ep_r)/(1-d)>r$, $H_M(M)=0$
and 
\[
0<\la u(\cdot;\mu,\ell)\ra<r\qquad (\ell>\ell_r),
\]
by taking $\mu_r=M$ in Lemma \ref{lem:avu}, 
we have a solution $\mu=\mu^*(\ell)\in[\omu+\ep_r, M]$ of the equation
$H_M(\mu)=\la u(\cdot;\mu,\ell)\ra$.
In the sequel, by putting $\ell_M=\ell_r$,  we obtain
the desired solution as $u^*(\cdot;\ell):=u(\cdot;\mu^*(\ell),\ell)$.
\eproof

\par\bigskip
\noindent
\underline{{\it Completion of the proof of Theorem \ref{thm:1}}}:
By \eqref{v} and $w=u$, we obtain the desired solution in
the theorem
by Lemma \ref{lem:Lr}.
As $\ell\to\infty$, $\la u^*(\cdot;\ell)\ra\to0$
follows from Lemma \ref{lem:ur}.
In view of
\[
\la u(\cdot;\mu^*(\ell),\ell)\ra =\frac{M-\mu^*(\ell)}{1-d},
\]
we assert 
$\lim_{\ell\to\infty}\mu^*(\ell)= M$.
This implies that the convergence of $u^*(\cdot;\ell)$ to
the homoclinic solution $u^h(\cdot)$ stated in Theorem \ref{thm:1},
which completes the proof of the theorem.
\eproof
\par\medskip
In the rest of this section we give the proof for Lemmas \ref{lem:rho}
and \ref{lem:avu} though the proof for Lemma \ref{lem:dxi} is
left in the later section.
\par\medskip
\noindent
\underline{{\it Proof of Lemma \ref{lem:rho}}}: 
We separate the integral of the right-hand side of \eqref{Lxi} as 
\begin{equation}
\label{intxi}
\begin{aligned}
&
\int_{\eta(\xi, \mu)}^{\xi}\frac{dz}{\sqrt{2(G(\eta(\xi,\mu); \mu)-G(z;\mu))}} \\
&
\quad =\int_{\eta(\xi, \mu)}^{\alpha(\mu)}\frac{dz}{\sqrt{2(G({\eta(\xi, \mu)}; \mu)-G(z;\mu))}} \\
&
\qquad +
\int_{\alpha(\mu)}^{\xi}\frac{dz}{\sqrt{2(G({\eta(\xi, \mu)}; \mu)-G(z;\mu))}}.
\end{aligned}
\end{equation}
Making use of
\begin{align*}
G(z;\mu)=G(z;\mu)-G(0;\mu)=\left(\int_0^1\int_0^1s {g}_u(\tau s z;\mu)ds d\tau\right)z^2< 0\quad(0<z<\gamma(\mu)), 
\end{align*}
where  ${g}_u(u;\mu)=\partial {g}(u;\mu)/\partial u$, and
\begin{align*}
\sqrt{G({\eta(\xi, \mu)};\mu)-G(z;\mu)}\le \sqrt{-G(z;\mu)}
=\sqrt{G(0;\mu)-G(z;\mu)}\quad(0<z<\gamma(\mu)),
\end{align*}
we have
\begin{align*}
\int_{\eta(\xi, \mu)}^{\alpha(\mu)}\frac{dz}{\sqrt{(G({\eta(\xi, \mu)}; \mu)-G(z;\mu))}}
&\ge \int_{\eta(\xi, \mu)}^{\alpha(\mu)}
\frac{dz}{\sqrt{(G(0;\mu)-G(z;\mu))}} \\
&
=\int_{\eta(\xi, \mu)}^{\alpha(\mu)}
\frac{dz}{z\sqrt{-\int_0^1\int_0^1 s {g}_u(\tau s z; \mu)ds d\tau}}. 
\end{align*}
Putting
\[
C_1(\mu_r):=\max_{\omu< \mu\leq \mu_r}\max_{0\le u\le \gamma(\mu)}|{g}_u(u; \mu)|,
\]
we have 
\begin{align*}
\int_{\eta(\xi, \mu)}^{\alpha(\mu)}\frac{dz}{z\sqrt{-\int_0^1\int_0^1s {g}_u(\tau s z; \mu)ds d\tau}} 
&\geq \frac1{\sqrt{2C_1(\mu_r)}}\log(\alpha(\mu)/{\eta(\xi, \mu)}),
\end{align*}
in the sequel, \eqref{intxi} implies
\begin{align*}
\int_{\eta(\xi, \mu)}^{\xi}\frac{dz}{\sqrt{2(G({\eta(\xi, \mu)}; \mu)-G(z;\mu))}}
>\frac1{\sqrt{2C_1(\mu_r)}}\log(\alpha(\omu)/{\eta(\xi, \mu)})
\end{align*}
follows
from the monotonicity of $\alpha(\mu)$ with respect to $\mu$.
In view of \eqref{ee} in the proof of Lemma \ref{lem:eta}
we have the desired assertion.
\eproof

\par\vspace{0.5cm}\noindent
\underline{{\it Proof of Lemma \ref{lem:avu}}}: 
We let $\xi^*=\xi(\mu,\ell)$ of Lemma \ref{lem:rell}.
The change of variable
$z=u(x; \mu,\ell)$ in the integral leads us to
\begin{align*}
\int_0^\ell ~u(x;\mu,\ell)\ dx=\int_{\xi^*}^{\eta(\xi^*, \mu)}  \frac{z dz}{u_x(x; \mu,\ell)} =
\int_{\eta(\xi^*, \mu)}^{\xi^*} \frac{zdz}{\sqrt{2(G({\eta(\xi^*, \mu)};\mu)-G(z;\mu))}}.
\end{align*}
Set
\[
\zeta_1:=\frac{1}{2}\min_{\mu\in[\omu+\ep_r, \mu_r]}
\min\{z\in [0, \alpha(\mu)): {g}_u(z;\mu)=0 \},
\]
and 
\begin{align}
\label{g1}
{g}_1:=\min_{(u, \mu)\in[0, \zeta_1]\times[\omu+\ep_r, \mu_r]}
|{g}_u(u;\mu)|>0.
\end{align}
We separate the integral as
\begin{align*}
&
\int_{\eta(\xi^*, \mu)}^{\xi^*}\frac{zdz}{\sqrt{2(G({\eta(\xi^*, \mu)}; \mu)-G(z;\mu))}} 
=I_1+I_2, \\
&
\qquad
I_1:=\int_{\eta(\xi^*, \mu)}^{\zeta_1}\frac{zdz}{\sqrt{2(G({\eta(\xi^*, \mu)}; \mu)-G(z;\mu))}},
\\
&
\qquad
I_2:=\int_{\zeta_1}^{\xi^*}\frac{zdz}{\sqrt{2(G({\eta(\xi^*, \mu)}; \mu)-G(z;\mu))}}.
\end{align*}
We simply write $\eta^*={\eta(\xi^*, \mu)}$. 
We first estimate $I_1$.
\begin{align*}
&
G(\eta^*;\mu)-G(z;\mu)
\\
&\qquad
=\left\{\int_0^1\left(\int_0^1
g_u(\tau(z+s(\eta^*-z));\mu)d\tau\right) (z+s(\eta^*-z))ds
\right\}(\eta^*-z)\\
&
\qquad
\geq
\min_{ (u, \mu)\in[\eta^*, z]\times [\omu+\ep_r, \mu_r]  }|{g}_u(u;\mu)|\left\{\int_0^1 (z+s(\eta^*-z))ds\right\}(\eta^*-z)\\
&
\qquad \geq
\frac{{g}_1}{2}(z+\eta^*)(z-\eta^*)
\geq \frac{{g}_1}{2}z(z-\eta^*),
\end{align*}
where ${g}_1$ is as in \eqref{g1}.
From this we see
\begin{align*}
I_1\leq
\int_{\eta^*}^{\zeta_1}\frac{z dz}{\sqrt{{g}_1z(z-\eta^*)}}
\leq \frac{1}{\sqrt{{g}_1}}\int_{\eta^*}^{\zeta_1}\frac{\sqrt{z}dz}{\sqrt{z-\eta^*}} 
\leq \frac{2\zeta_1}{\sqrt{{g}_1}},
\end{align*}
which is uniformly bounded in $[\omu+\ep_r, \mu_r]\times[\ell_0,\infty)$.

As for $I_2$, in view of
\begin{align*}
G(\xi^*; \mu)-G(z;\mu)=\left(\int_0^1 g(z+s(\xi^*-z);\mu)ds\right)(\xi^*-z),
\end{align*}
and the change of variable
$z=\zeta_1+(\xi^*-\zeta_1)\sin^2\theta$,
we see
\begin{align*}
&
\int_{\zeta_1}^{\xi^*}\frac{zdz}{\sqrt{2(G(\xi^*; \mu)-G(z;\mu))}}\\
&
=\int_0^{\pi/2}\frac{\sqrt2\{\zeta_1+(\xi^*-{\zeta_1})\sin^2\theta\}(\xi^*-\zeta_1)\sin\theta\cos\theta d\theta}{\sqrt{\left(\int_0^1 {g}(\zeta_1+(\xi^*-\zeta_1)\sin^2\theta+s(\xi^*-\zeta_1)\cos^2\theta;\mu )ds\right)(\xi^*-\zeta_1)\cos^2\theta}} \\
&
=\int_0^{\pi/2}\frac{\sqrt2\{{\zeta_1}+(\xi^*-{\zeta_1})\sin^2\theta\}\sqrt{\xi^*-\zeta_1}\sin\theta d\theta}{\sqrt{\int_0^1 {g}(\zeta_1+(\xi^*-\zeta_1)\sin^2\theta+s(\xi^*-\zeta_1)\cos^2\theta;\mu )ds}}.
\end{align*}
Since $\xi^*=\xi(\mu,\ell)$, $\I_2$ is also uniformly bounded in 
$[\omu+\ep_r, \mu_r]\times[\ell_0,\infty)$.
We proved Lemma \ref{lem:avu}.
\eproof

\par\bigskip
\section{Proof of (i) and (ii) of Lemma \ref{lem:1}}\label{proof_lemma}
In the previous section we constructed the monotone solution in the parameter regime $\mu>\omu$. 
On the other hand, noticing $G(\alpha(\mu_c))=G(\beta(\mu_c))<0$, there is
$\omu<\mu_1$ such that
\begin{align*}
G(\beta(\omu))=0>G(\beta(\mu))\quad (\mu\in(\mu_c, \omu)),
\end{align*}
where $\mu_1$ is the number defined in \eqref{mu1}.

As for $\mu\in(\mu_c, \omu]$, we can also construct
a monotone solution to \eqref{ugx}.
Indeed, under this regime for $\mu$ there is a homoclinic solution joining $\beta(\mu)$ to itself for $\mu\in(\mu_c, \omu)$ or a heteroclinic solutions joining $0(=\alpha(\omu))$ and $\beta(\omu)$ for $\mu=\omu$ to the equation
\[
du_{xx}+{g}(u;\mu)=0 \quad (-\infty<x<\infty),
\]
thus it is expected that there is a solution to \eqref{ugx} close to the homoclinic/heteroclinic solution for sufficiently large $\ell$.


\begin{rmk}
{\em
Unlike the arguments in \S\ref{proof_theorem}, we don't discuss  a similar uniform dependence on $\mu$ as shown for $\xi(\mu,\ell)$ in \eqref{rell},
since we do not solve the equation \eqref{Mmuu} for given $M$. 
Indeed, this problem is much more difficult. We need a further study in this case.}
\end{rmk}

We let $\omega_*(\mu)$ be a number satisfying
\[
G(\omega_*(\mu); \mu)=G(\beta(\mu);\mu), \quad \omega_*(\mu)\in(0, \alpha(\mu)),
\]
and set
\begin{align*}
u(0)=\omega\in(\omega_*(\mu), \alpha(\mu)).\qquad 
\end{align*}

\begin{figure}[bt]
\begin{tabular}{lr}
\begin{minipage}[t]{0.5\hsize}
\centering
\includegraphics[keepaspectratio, scale=0.4]{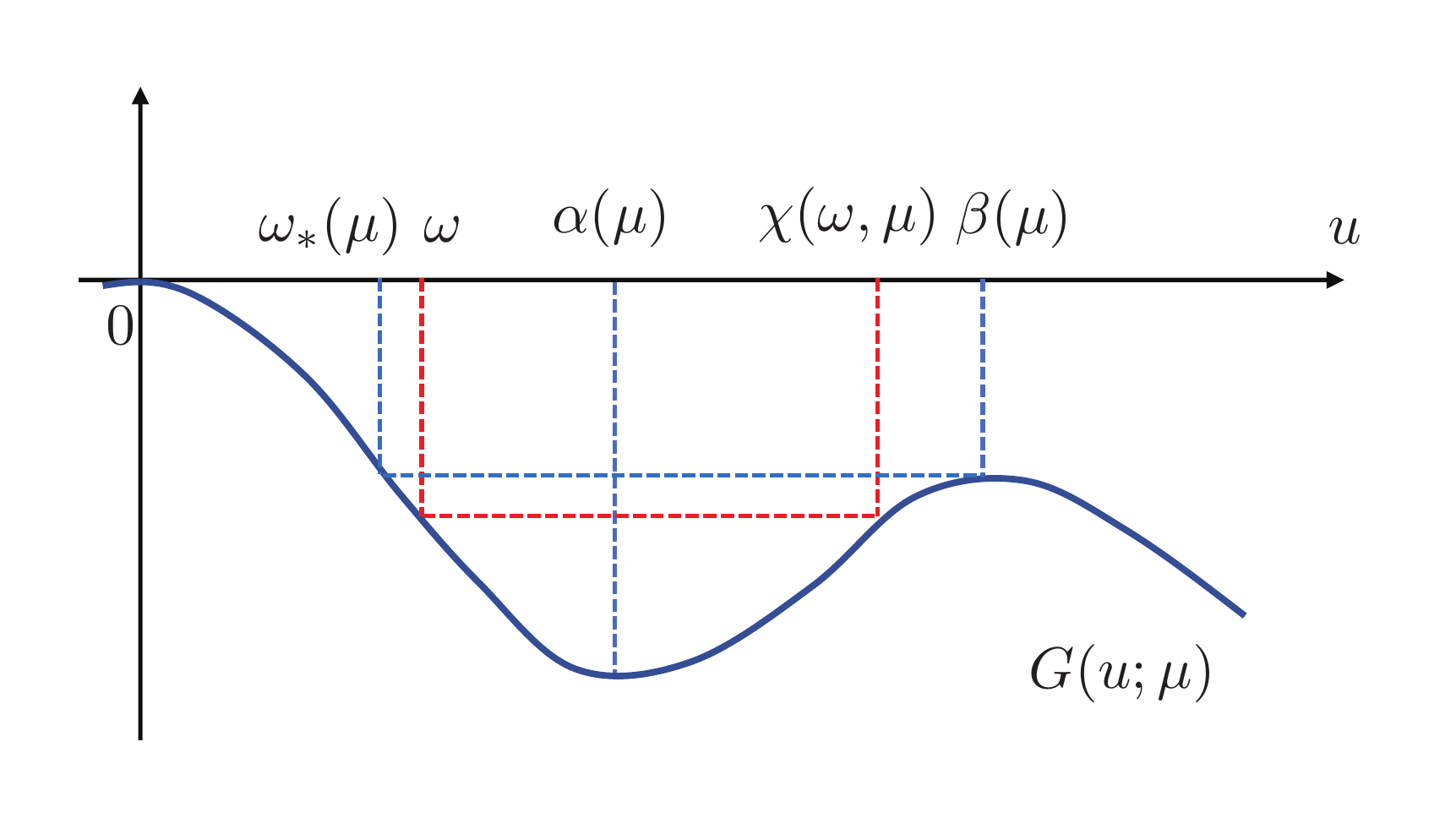}
\caption{\small{Graph of $G(u;\mu)$}}
\label{fig5}
\end{minipage}
&
\begin{minipage}[t]{0.5\hsize}
\centering
\includegraphics[keepaspectratio, scale=0.45]{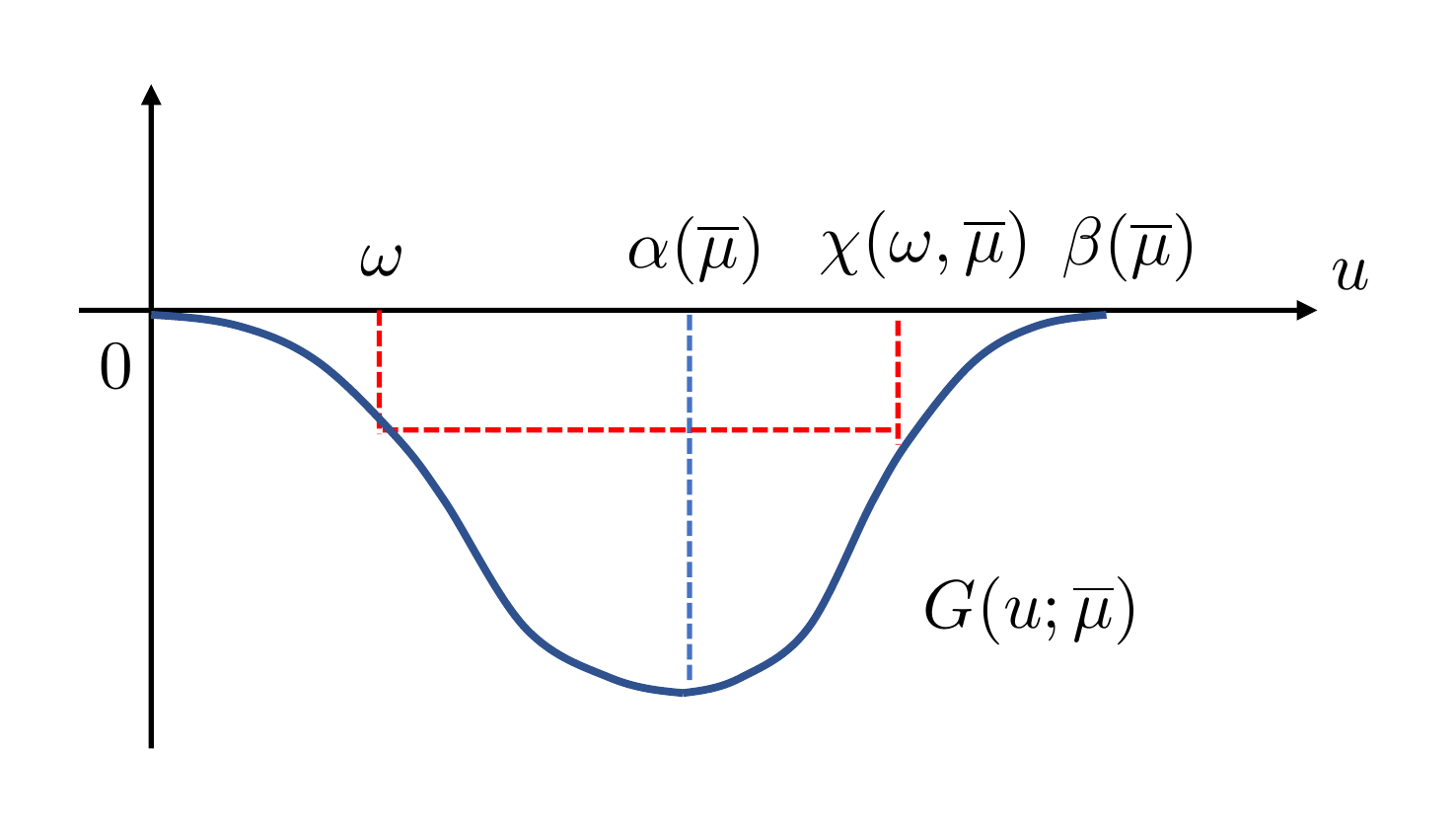}
\caption{\small{Graph of $G(u;\omu)$}}
\label{fig6}
\end{minipage}
\end{tabular}
\end{figure}

Let $\chi=\chi(\omega, \mu)$ be the value solving the equation
\[
G(\chi;\mu)=G(\omega; \mu)
\]
(see Figure \ref{fig5}).
We note that $\omega_*(\omu)=0$ since $G(\beta(\omu);\omu)=0$ 
and $\chi(\omega_*(\mu),\mu)=\beta(\mu)$
(see  Figure \ref{fig6}).

\begin{figure}[bt]
\centering
\includegraphics[keepaspectratio, scale=0.6]{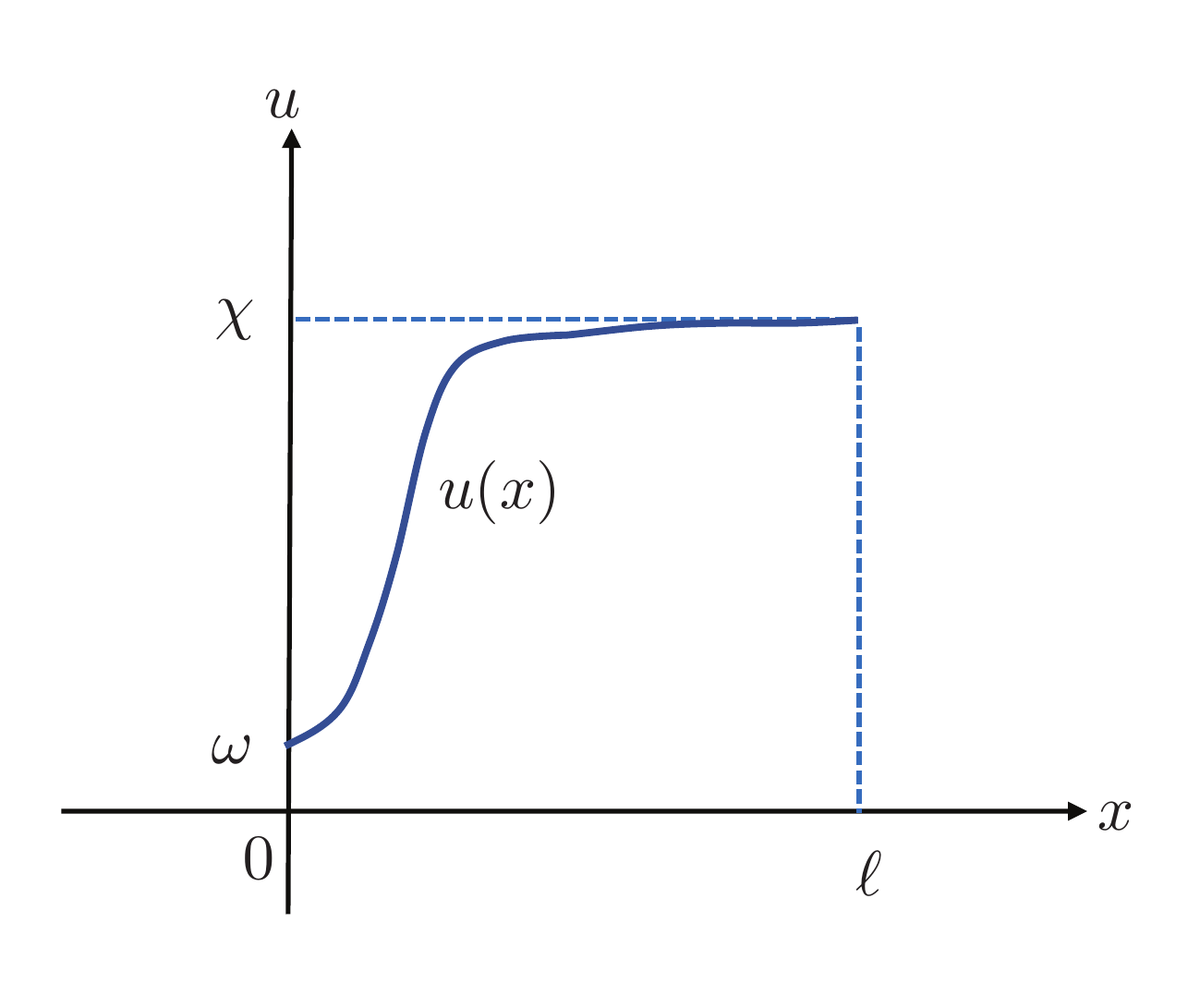}
\caption{\small{Profile of the solution $u(x)$}}
\label{fig7}
\end{figure}

Given $\mu\in(\mu_c, \omu]$, we look for a monotone increasing solution satisfying
$u(\ell)=\chi (\omega,\mu)$ (see Figure \ref{fig7}).

In a similar way to obtain \eqref{Lxi} we reduce the problem to finding a solution of 
\begin{equation}
\label{Leta}
\begin{aligned}
\ell=
\tilde\rho(\omega, \mu):=&\int_{\omega}^{\chi(\omega,\mu)}
\frac{dz}{\sqrt{2(G(\omega; \mu)-G(z;\mu))}}\\
=
&\int_{\omega}^{\chi(\omega, \mu)}
\frac{dz}{\sqrt{2(G(\chi(\omega,\mu); \mu)-G(z;\mu))}}.
\end{aligned}
\end{equation}
In order to show the existence of a solution to \eqref{Leta}, we first show
\par\medskip\noindent
\begin{lem}
\label{lem:trho}
Given $\mu\in(\mu_c, \omu]$, $\lim_{\omega\to \omega_*(\mu)}\tilde\rho(\omega, \mu)=\infty$ holds.
\end{lem}
\proof
Put
\[
C_1:=\max_{{{\mu_c}}\le \mu\leq\omu}\max_{\alpha(\mu)\le u\le \beta(\mu)}|{g}_u(u; \mu)|,
\]
and we obtain 
\begin{align*}
&
G(\chi(\omega,\mu); \mu)-G(z; \mu)
\leq G(\beta(\mu); \mu)-G(z; \mu)\\
&\quad=\left(\int_0^1[-{g}(\beta(\mu)+s(z-\beta(\mu)); \mu)]ds\right)(z-\beta(\mu))\\
&\quad=\left(\int_0^1[{g}(\beta(\mu)+s(z-\beta(\mu));\mu)-{g}(\beta(\mu);\mu)]ds\right)(\beta(\mu)-z) \\
&\quad=-\left(\int_0^1\int_0^1
sg_u(\beta(\mu)+\tau s(z-\beta(\mu));\mu)\ ds d\tau\right)(\beta(\mu)-z)^2\\
&\quad\leq C_1(\beta(\mu)-z)^2
\end{align*}
for $z \in [\alpha(\mu), \chi(\omega, \mu)]$.
Making use of the above inequality, we have
\begin{equation}
\label{tilrho}
\begin{aligned}
\tilde\rho(\omega, \mu)&\geq 
\int_{\alpha(\mu)}^{\chi (\omega,\mu)}
\frac{dz}{\sqrt{2(G(\chi (\omega,\mu); \mu)-G(z;\mu))}}\\
&
\geq
\int_{\alpha(\mu)}^{\chi (\omega,\mu)}\frac{dz}{\sqrt{2C_1(\beta(\mu)-z)^2}}
=\frac1{\sqrt{2C_1}}\log\frac{\beta(\mu)-\alpha(\mu)}{\beta(\mu)-\chi(\omega,\mu)}.
\end{aligned}
\end{equation}
Since $\lim_{\omega\to \omega_*(\mu)}\chi (\omega,\mu)=\beta(\mu)$, we obtain the desired assertion.
\eproof
\par\medskip

\begin{lem}
\label{lem:dtrho}
For each $\mu\in(\mu_c, \omu]$
there exist $\delta_1>0$ such that if $\omega\in (\omega_*(\mu), \omega_*(\mu)+\delta_1]$, then
\[
\frac{\partial\tilde\rho}{\partial\omega}(\omega,\mu)<0
\]
holds.
\end{lem}



\par\medskip\noindent
The proof is given in \S\ref{sec:proof}.

Applying 
Lemmas \ref{lem:trho} and \ref{lem:dtrho}, 
we have the next lemma.

\begin{lem}
\label{lem:trhosol}
For each $\mu\in(\mu_c, \omu]$ consider the equation
$\tilde\rho(\omega, \mu)=\ell$. Then
there exists $\tilde\ell>0$ such that for $\ell>\tilde\ell$ the equation allows a unique solution $\omega=\omega(\mu, \ell)$ 
satisfying $\lim_{\ell\to\infty}\omega(\mu, \ell)=\omega_*(\mu)$.
\end{lem}


\par\bigskip\noindent
{\it Proof of \rm{(i)} and \rm{(ii)} in Lemma \ref{lem:1}}:
The solution (i) in Lemma \ref{lem:1} is given by the solution to 
$\tilde\rho(\omega, \mu)=\ell$ of \eqref{Leta}, namely, the solution 
to 
\[
du_{xx}+g(u;\mu)=0~~(x\in(0,\ell)), \qquad u(0)=\omega(\mu,\ell)
\]
is nothing but $u(\cdot; \mu, \ell)$ in (i).  
Noticing that the homoclinic solution to \eqref{ugb} is obtain
by the solution to 
\[
du_{xx}+g(u;\mu)=0~~(x\in\bR), \qquad u(0)=\omega_*,
\]
we have the assertion of (i) on the convergence of the solution $u(x;\mu, \ell)$ 
to the homoclinic solution as $\ell\to\infty$ is valid by Lemma \ref{lem:trhosol}.

As for the case $\mu=\omu$, we similarly obtain the solution $u(\cdot;\omu, \ell)$
in (ii) by  one of 
\[
du_{xx}+g(u;\omu)=0~~(x\in(0,\ell)), \qquad u(0)=\omega(\omu, \ell).
\]
We let the solution shift as $\tilde{u}(x;\omu,\ell)=u(x+\ell/2; \omu, \ell)$. Then
\begin{align*}
& 
d\tilde{u}_{xx}+g(\tilde{u};\omu)=0~~(x\in(-\ell/2,\ell/2)), \\
&
\tilde{u}(-\ell/2; \omu, \ell)=\omega(\omu, \ell), \qquad 
\tilde{u}(\ell/2; \omu, \ell)=\chi(\omega(\omu, \ell),\omu).
\end{align*}
Since 
\[
\lim_{\ell\to\infty}\omega(\omu, \ell)=0,\qquad \lim_{\ell\to\infty}\chi(\omega(\omu, \ell),\omu)=\beta(\omu),
\]
we have the assertion for the convergence to the heteroclinic solution to \eqref{sol:hetero}.
This concludes the proof.
\eproof
\par\medskip
In order to investigate the profile the non-constant stationary solution numerically,
we calculate the following auxiliary evolutional system 
which is reduced from the original system \eqref{eq:NSI}:
\begin{equation}\label{eq:auxi_NS}
\left\{
\begin{aligned}
&
N_t=\frac{N^2S}{1+(k_N/k_I) N}-N+D_N\Delta N, \\
&
S_t=-\frac{N^2S}{1+(k_N/k_I) N}+N+\Delta S.  \\
\end{aligned}
\right.
\end{equation}
The initial conditions are given as follows:
\begin{align*}
&N(0,x)=\left\{
    \begin{aligned}
    &2A &  \Big(x \in \Big[\frac{L}{4}, \frac{3L}{4}\Big)\Big), \\
    &0  & (\text{otherwise}).  
    \end{aligned}
    \right.\\
&S(0,x)=0 \quad (x \in(0, L)).
\end{align*}
We note that any stationary solution to \eqref{eq:NSI} with $D_I=0$ is obtained by a stationary solution $(N^*(x), S^*(x))$ of the above system as put $(N,S,I)=(N^*(x), S^*(x), N^*(x))$.

\begin{figure}[bt]
	\begin{center}
	\begin{tabular}{cc}
		\includegraphics[width=5cm, bb=0 0 400 400]{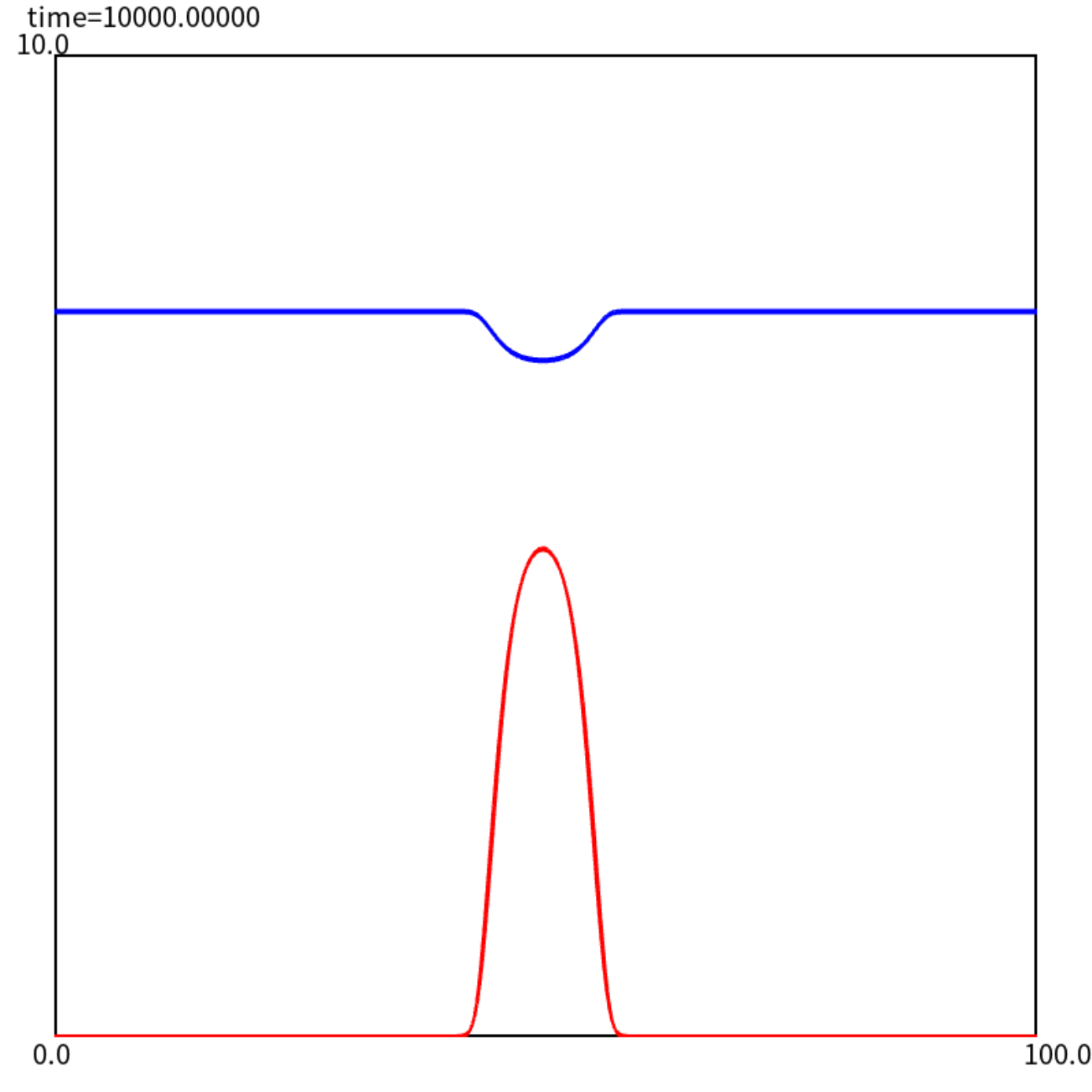}
		&\includegraphics[width=5cm, bb=0 0 400 400]{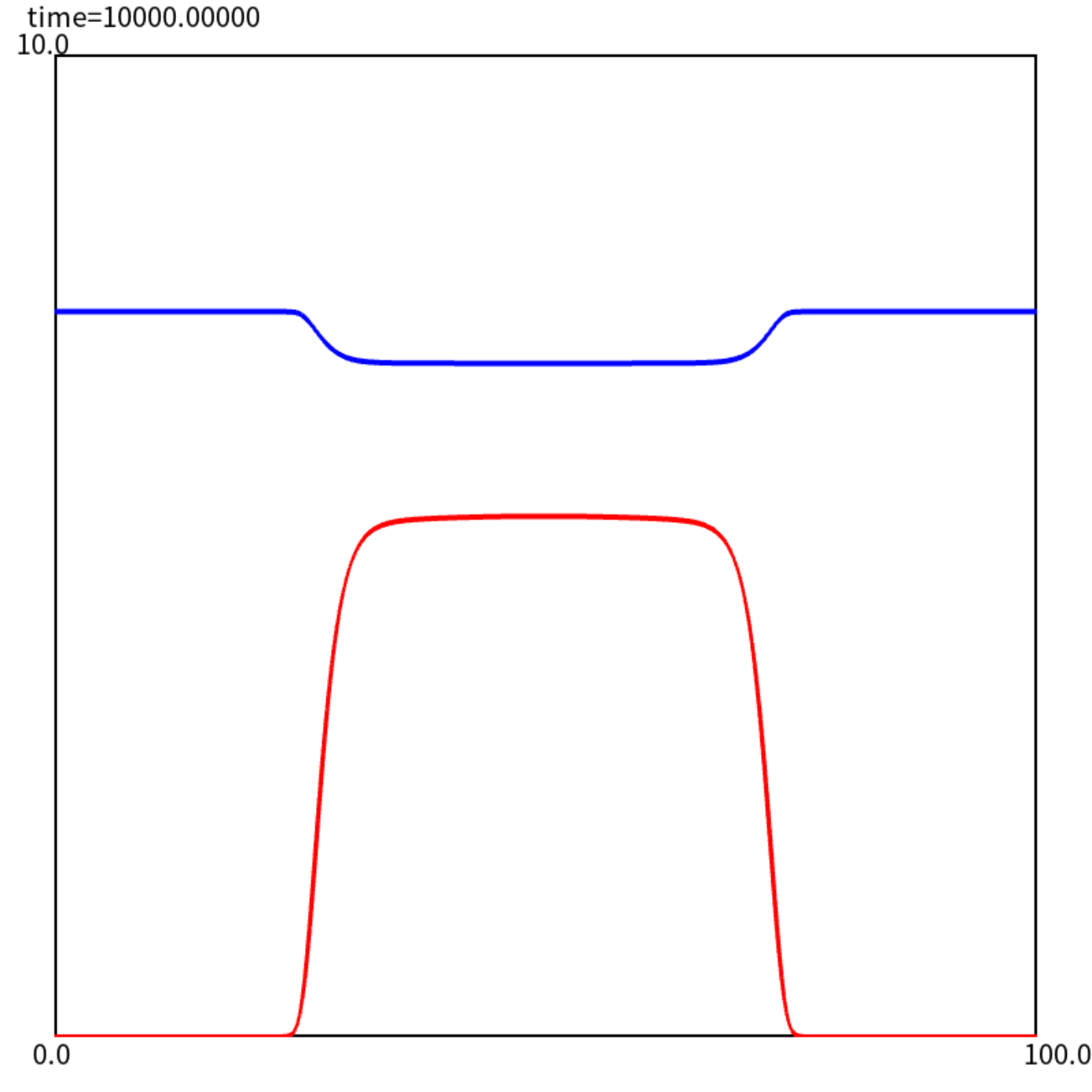}\\
		(a) & (b)
	\end{tabular}
	\end{center}
	\caption{\small{The profiles of the non-constant stationary solution given by a numerical simulation for 
	             the parameters, $k_N=2.0, k_I=0.3, D_N=0.1, K=20/3, A=9.5$.  
	             We compute \eqref{eq:auxi_NS} until $t=10000$ for a suitable initial condition.
			The red and blue curves correspond to the solution of $N$ and $S$ to \eqref{eq:auxi_NS} in the interval $[0, 100]$, respectively. 
	            (a) $A=7.8$ and (b) $A=9.5$.}
}
	\label{fig8}
\end{figure}

We calculate the solution of evolutional problem \eqref{eq:auxi_NS} until $t=10000$ to obtain Figure \ref{fig8}.
As seen in the figure, according to the value of $A$, the system numerically exhibits two types of stationary solutions which correspond to the ones in (i) and (iii) of Lemma \ref{lem:1}.  Furthermore, we remark that in our parameter setting, the corresponding stationary solutions in the three-component system \eqref{eq:NSI} with $D_I=0$ are destabilized, though not displayed here.

\noindent

\section{Proof of Lemmas \ref{lem:dxi} and \ref{lem:dtrho}}\label{appendix}\label{sec:proof}
\subsection{Proof of Lemma \ref{lem:dxi}}
First consider $\partial\eta/\partial\xi$.
In view of \eqref{GG} we have
\begin{align}
\label{etag}
\frac{\partial\eta}{\partial\xi}(\xi,\mu)=\frac{ g(\xi; \mu)}{ g(\eta(\xi,\mu);\mu)}. \qquad
\end{align}
Applying \eqref{ee} to the denominator of the right-hand side yields that
arbitrarily given $N>0$, there is $\delta_N>0$, independent of $\mu\in (\omu, \mu_r]$, such that for $\mu\in (\omu, \mu_r]$
\begin{align}
\label{eN}
\frac{\partial\eta}{\partial\xi}(\xi,\mu)<-N \quad \text{holds}\quad \text{if} \quad 
\xi\in(\gamma(\mu)-\delta_N, \gamma(\mu)).
\end{align}

\par\medskip
Go to the next step. 
Put
\[
\zeta:=\min_{\mu\in(\omu, \mu_r]}
\min\{z\in [0, \alpha(\mu)): g_u(z;\mu)=0 \},
\]
and take $\xi_m$ so that $\eta(\xi_m,\mu)<\zeta$.
Then
\begin{align}
\label{xim}
0>g(\eta(\xi,\mu);\mu)>g(\zeta;\mu) \qquad(\xi\in[\xi_m, \gamma(\mu))).
\end{align}
For $\xi\in[\xi_m, \gamma(\mu))$, 
we divide the right-hand side of \eqref{Lxi} as
%
\begin{align*}
\rho(\xi,\mu)&=\frac1{\sqrt2}(\I_1+\I_2),  \\
\I_1&:=\int_{\eta(\xi, \mu)}^{\zeta}
\frac{dz}{\sqrt{G(\eta(\xi, \mu); \mu)-G(z;\mu)}} \\
&=\int_0^1
\frac{(\zeta-\eta(\xi,\mu)) d\tau}{\sqrt{G(\eta (\xi,\mu);\mu)
- G(\eta (\xi,\mu)+\tau(\zeta-\eta (\xi,\mu));\mu)} },\\
\I_2&:=
\int_{\zeta}^{\xi}\frac{dz}{\sqrt{G({\eta (\xi,\mu)}; \mu)-G(z;\mu)}}
=\int_{\zeta}^{\xi}\frac{dz}{\sqrt{G(\xi; \mu)-G(z;\mu))}}.
\end{align*}
In the computation below
we use the following formula for $\Psi(y,z)$ 
satisfying $\lim_{\tau\to0}\tau\Psi(y, y+\tau(a-y))=0$:
\begin{equation}
\label{Phi}
\begin{aligned}
&
\frac{\partial}{\partial y}\int_0^1 (a-y)\Psi(y, y+\tau(a-y))d\tau 
=-\Psi(y,a)\\
&
+\int_0^1(a-y)\left[\frac{\partial}{\partial y}\Psi(y, y+\tau(a-y))+\frac{\partial}{\partial z}\Psi(y, y+\tau(a-y))\right]d\tau,
\end{aligned}
\end{equation}
where 
\[
\Psi(y, y+\tau(a-y)), \quad
\frac{\partial}{\partial y}\Psi(y, y+\tau(a-y)), \quad
\frac{\partial}{\partial z}\Psi(y, y+\tau(a-y))
\]
are integrable in $\tau\in(0,1)$.
Indeed, 
\begin{align*}
&
\frac{\partial}{\partial y}\int_0^1(a-y)\Psi(y, y+\tau(a-y))d\tau
=-\int_0^1\Psi(y, y+\tau(a-y)) d\tau \\
&
+
\int_0^1(a-y)\left[\frac{\partial}{\partial y}\Psi(y, y+\tau(a-y))+(1-\tau)\frac{\partial}{\partial z}\Psi(y, y+\tau(a-y))\right]d\tau.
\end{align*}
Since
\begin{align*}
&
-\int_0^1(a-y)\tau\frac{\partial}{\partial z}\Psi(y, y+\tau(a-y))d\tau =-
\int_0^1\tau\frac{\partial}{\partial\tau}\Psi(y, y+\tau(a-y))d\tau \\
&
=-[\tau\Psi(y, y+\tau(a-y))]_{\tau=0}^{1} +\int_0^1\Psi(y, y+\tau(a-y))d\tau,
\end{align*}
it is easy to see that the identity of \eqref{Phi} is true.

Putting
\begin{align*}
\Psi(y,z)=\frac{1}{\sqrt{G(y; \mu)-G(z;\mu)}}, 
\end{align*}
we compute
\begin{align*}
\frac{\partial \I_1}{\partial\xi}
 =\frac{\partial}{\partial\xi}\int_0^1
 (\zeta-\eta(\xi,\mu))\Psi(\eta(\xi,\mu), \eta(\xi,\mu)+\tau(\zeta-\eta(\xi,\mu))) d\tau.
 \end{align*}
By \eqref{Phi} (with $a=\zeta$, $y=\eta(\xi, \mu)$,  and $z= \eta(\xi,\mu)+\tau(\zeta-\eta(\xi,\mu))$), 
\begin{align*}
\frac{\partial\Psi}{\partial y}(y,z)=-\frac{ g(y;\mu)}{2[G(y; \mu)-G(z;\mu)]^{3/2}}, \quad
\frac{\partial\Psi}{\partial z}(y,z)=\frac{ g(z;\mu)}{2[G(y; \mu)-G(z;\mu)]^{3/2}}, 
\end{align*}
\[
\lim_{\tau\to+0}\frac{\tau}{\sqrt{G(\eta(\xi,\mu);\mu)- G(\eta(\xi,\mu)+
\tau(\zeta-\eta(\xi,\mu));\mu)} }=0\quad
(\text{from L'H\^{o}pital's rule}),
\]
and the chain rule of differentiation, 
we have
\begin{equation}
\label{I1}
\begin{aligned}
&
 \frac{\partial \I_1}{\partial\xi}
 =-\frac{\partial\eta}{\partial\xi}
 \left\{  
 \frac{1}{\sqrt{G(\eta(\xi,\mu);\mu)- G(\zeta;\mu)} } 
 \right. \\
 &
 \quad
 \left.
 +\int_0^1\frac{(\zeta-\eta(\xi,\mu))\{ g(\eta(\xi,\mu);\mu) - g(\eta(\xi,\mu)+\tau(\zeta-\eta(\xi,\mu));\mu)\}}{2[G(\eta(\xi,\mu);\mu)- G(\eta(\xi,\mu)+\tau(\zeta-\eta(\xi,\mu));\mu)]^{3/2} } d\tau
 \right\}
\\
&\quad
>  \frac{-\partial\eta/\partial\xi}{\min_{\mu\in(\omu, \mu_r]}\sqrt{-G(\zeta;\mu)} }>0,
\end{aligned}
\end{equation}
where we used \eqref{xim}.
We note the integral of the above equality
is positive by \eqref{etag}. 

We next handle $\I_2$.
Write
\begin{align*}
\I_2&=\J_1+\J_2, \\
\J_1&:=\int_\zeta^{\alpha(\mu)}  \frac{1}{\sqrt{G(\xi;\mu)- G(z;\mu)} } dz, 
\qquad
\J_2:=\int_{\alpha(\mu)}^\xi  \frac{1}{\sqrt{G(\xi;\mu)- G(z;\mu)} } dz.
\end{align*}
We have
\begin{align*}
\frac{\partial\J_1}{\partial\xi}=
\frac1{2}\int_\zeta^{\alpha(\mu)} \frac{ g(\xi;\mu)}{[G(\xi;\mu)-G(z;\mu)]^{3/2}}dz<0.
\end{align*}
We estimate the bound of this one as
\begin{align*}
\left|\frac{\partial\J_1}{\partial\xi}\right|\leq 
\frac{| g(\xi;\mu)|}{2}\int_\zeta^{\alpha(\mu)} \frac{1}{\{G(\xi;\mu)-G(\zeta;\mu)\}^{3/2}}dz
=\frac{| g(\xi;\mu)|(\alpha(\mu)-\zeta)}{2\{G(\xi;\mu)-G(\zeta;\mu)\}^{3/2}},
\end{align*}
hence,
\begin{align*}
\left|\frac{\partial\J_1}{\partial\xi}\right|\leq 
c_1:=
\frac{\max_{(u,\mu)\in[\zeta,  \alpha(\bar \mu)]\times(\omu, \mu_r]}| g(u;\mu)|( \alpha(\bar \mu)-\zeta)}{2\min_{\mu\in(\omu, \mu_r]}\{G(\xi_m; \mu)-G(\zeta;\mu)\}^{3/2}},
\end{align*}
which yields
\begin{align}
\label{J1}
\frac{\partial\J_1}{\partial\xi}>-c_1.
\end{align}

On the other hand, we write 
\begin{align*}
\J_2=\int_{\alpha(\mu)}^\xi  \frac{1}{\sqrt{G(\xi;\mu)- G(z;\mu)} } dz
=\int_0^1\frac{\xi-\alpha(\mu)}{\sqrt{G(\xi;\mu)- G(\xi-\tau(\xi-\alpha(\mu)); \mu)} } d\tau,
\end{align*}
and compute
\begin{align*}
\frac{\partial\J_2}{\partial\xi}&=
\int_0^1\frac{1}{\sqrt{G(\xi;\mu)- G(\xi-\tau(\xi-\alpha(\mu)); \mu)} } d\tau \\
&
\quad
-\frac{\xi-\alpha(\mu)}{2}\int_0^1 
\frac{g(\xi; \mu)-g(\xi-\tau(\xi-\alpha(\mu));\mu)(1-\tau)}{[G(\xi;\mu)- G(\xi-\tau(\xi-\alpha(\mu)); \mu)]^{3/2}} d\tau\\
&
=\int_0^1\frac{1}{\sqrt{G(\xi;\mu)- G(\xi-\tau(\xi-\alpha(\mu)); \mu)} } d\tau  \\
&\quad
-\frac{\xi-\alpha(\mu)}{2}\int_0^1\frac{g(\xi; \mu)-g(\xi-\tau(\xi-\alpha(\mu));\mu)}{[G(\xi;\mu)- G(\xi-\tau(\xi-\alpha(\mu)); \mu)]^{3/2}} d\tau \\
&\quad
-\frac{\xi-\alpha(\mu)}{2}\int_0^1 
\frac{\tau g(\xi-\tau(\xi-\alpha(\mu));\mu)}{[G(\xi;\mu)- G(\xi-\tau(\xi-\alpha(\mu)); \mu)]^{3/2}} d\tau.
\end{align*}
We estimate
\begin{align*}
&
\frac{| g(\xi; \mu) - g(\xi-\tau(\xi-\alpha(\mu));\mu)|}{[G(\xi; \mu)- G(\xi-\tau(\xi-\alpha(\mu)); \mu)]^{3/2}} \\
&
=\frac{\left|\int_0^1g_u(\xi-(1-s)(\xi-\alpha(\mu))\tau; \mu) ds\right|
(\xi-\alpha(\mu))\tau}{
\left[\left(
\int_0^1 g(\xi-(1-s)(\xi-\alpha(\mu))\tau; \mu) ds
\right)(\xi-\alpha(\mu))\tau\right]^{3/2}
}.
\end{align*}
We set
\begin{align*}
&
c_2:=\max_{\mu\in(\omu, \mu_r]}\max_{0\leq u\leq\gamma(\mu)}| g_u(u;\mu)|, 
\quad
c_3:=\min_{\mu\in(\omu,\mu_r]}
\min_{\frac{\alpha(\mu)+\xi_m}{2}\leq u\leq \gamma(\mu)} g(u;\mu), \\
&
c_4 :=\max_{\mu\in(\omu, \mu_r]}\max_{\alpha(\mu)\leq u\leq\gamma(\mu)}| g(u;\mu)|.
\end{align*}
Then 
\[
| g_u(\xi-(1-s)(\xi-\alpha(\mu))\tau; \mu)|\leq c_2. 
\]
On the other hand, since  
\[
\frac{\alpha(\mu)+\xi}{2}
\leq\xi-(1-s)(\xi-\alpha(\mu))\tau\leq\xi
\]
holds for $s\in[1/2,1]$, 
we have the estimate
\begin{align*}
&
\left|\int_0^1g(\xi-(1-s)(\xi-\alpha(\mu))\tau; \mu) ds\right|\geq
\int_{1/2}^1\min_{\frac{\alpha(\mu)+\xi_m}{2}\leq u\leq\gamma(\mu) } g(u;\mu)ds
\geq
\frac1{2}c_3.
\end{align*}
Thus,
\begin{equation}
\label{J22}
\begin{aligned}
&
\frac{\xi-\alpha(\mu)}{2}\int_0^1\frac{| g(\xi; \mu) - g(\xi-\tau(\xi-\alpha(\mu));\mu)|}{[G(\xi; \mu)- G(\xi-\tau(\xi-\alpha(\mu)); \mu)]^{3/2}}d\tau\\
&
\leq
\frac{\xi-\alpha(\mu)}{2}
\left(\int_0^1\frac{c_2(\xi-\alpha(\mu))\tau}{[
(c_3/2)(\xi-\alpha(\mu))\tau]^{3/2}})d\tau\right) \\
&
=\frac{c_2\sqrt2\sqrt{\xi-\alpha(\mu)}}{c_3^{3/2}}\left(\int_0^1\frac{d\tau}{\sqrt\tau} \right) 
=\frac{2\sqrt2 c_2\sqrt{\xi-\alpha(\mu)}}{c_3\sqrt{c_3}}.\end{aligned}
\end{equation}
\par\medskip
Finally, we estimate
\begin{equation}
\label{J23}
\begin{aligned}
&
\frac{\xi-\alpha(\mu)}{2}\int_0^1 
\frac{\tau g(\xi-\tau(\xi-\alpha(\mu));\mu)}{[G(\xi; \mu)- G(\xi-\tau(\xi-\alpha(\mu)); \mu)]^{3/2}} d\tau \\
&
\le
\frac{\xi-\alpha(\mu)}{2}
\int_0^1 
\frac{c_4\tau}{[(c_3/2)(\xi-\alpha(\mu))\tau]^{3/2}} d\tau
=\frac{\sqrt2 c_4\sqrt{\xi-\alpha(\mu)}}{c_3\sqrt{c_3}}.
\end{aligned}
\end{equation}

Combining \eqref{I1}, \eqref{J1}, \eqref{J22} and \eqref{J23}
yields
\begin{align*}
\frac{\partial\rho}{\partial\xi}>
-\frac1{c_0}\frac{\partial\eta}{\partial\xi}
-c_1-\tilde{c}_2-\tilde{c}_3,
\end{align*}
where we put
\begin{align*}
c_0:=\min_{\mu\in(\omu, \mu_r]}\sqrt{-G(\zeta;\mu)}, \quad
\tilde{c}_2:=
\frac{2\sqrt2 c_2\sqrt{\xi-\alpha(\omu)}}{c_3\sqrt{c_3}},
\quad 
\tilde{c}_3:=
\frac{2\sqrt2 c_4\sqrt{\xi-\alpha(\omu)}}{c_3\sqrt{c_3}}.
\end{align*}
We apply \eqref{eN}. Set
\[
N:=2c_0c, 
\qquad
c:=c_1+\tilde{c}_2+\tilde{c}_3, 
\]
and there is $\delta_N>0$ such that
for $\mu\in(\omu, \mu_r]$,
\[
\frac{\partial\rho}{\partial\xi}>
\frac1{c_0}N-c_1-\tilde{c}_2-\tilde{c}_3=c
\]
holds if $\xi\in(\gamma(\mu)-\delta_N, \gamma(\mu))\cup(\xi_m, \gamma(\mu))$.
Hence, we can take $\xi_0(\mu):=\max\{\gamma(\mu)-\delta_N, \xi_m\}$
so that \eqref{rcM} holds. 
This completes the proof of Lemma \ref{lem:dxi}.
\eproof

\par\vspace{1cm}\noindent

\subsection{Proof of Lemma \ref{lem:dtrho}}
Since a similar but simpler computations to those in the previous case does work, 
we only give a sketch of the proof. 
Put
\begin{align*}
\tilde\zeta:=
\max\{z_0: g_u(z;\mu)<0 ~(z_0<z\leq \beta(\mu)\},
\end{align*}
and verify
\begin{align*}
\frac{\partial \chi}{\partial\omega}=
\frac{ g(\omega;\mu)}{ g(\chi(\omega,\mu); \mu)}<0, \qquad
\lim_{\omega\to \omega_*(\mu)+0}\frac{\partial \chi}{\partial\omega}=-\infty.
\end{align*}
We separate the integral \eqref{Leta} as
\begin{align*}
\tilde\rho(\omega;\mu)&=\frac1{\sqrt2}(\tilde\I_1+\tilde\I_2),  \\
\tilde\I_1&:=\int_{\omega}^{\tilde\zeta}
\frac{dz}{\sqrt{G(\omega; \mu)-G(z;\mu)}} \\
\tilde\I_2&:=
\int_{\tilde\zeta}^{\chi (\omega,\mu)}\frac{dz}{\sqrt{G( \chi (\omega,\mu); \mu)-G(z;\mu)}} \\
&=\int_0^1
\frac{(\chi(\omega, \mu)-\tilde\zeta) d\tau}{\sqrt{G(\chi(\omega, \mu);\mu)
- G(\chi(\omega, \mu)+\tau(\tilde\zeta-\chi(\omega, \mu));\mu)} }.
\end{align*}
Using \eqref{Phi} and considering the definition of $\tilde\zeta$ and $\partial\chi/\partial\omega<0$, 
we compute
\begin{align*}
\frac{\partial\tilde\J_2}{\partial\omega}
&
=\frac{\partial\chi}{\partial\omega}\left\{\frac{1}{\sqrt{G(\chi;\mu)- G(\chi-\tau(\chi-\tilde\zeta); \mu)} }\right.  \\
&\qquad \qquad
\left.
-\frac{\chi-\tilde\zeta}{2}\int_0^1\frac{ g(\chi; \mu) - g(\chi-\tau(\chi-\tilde\zeta);\mu)}{[G(\chi;\mu)- G(\chi-\tau(\chi-\tilde\zeta); \mu)]^{3/2}} d\tau\right\} \\
&\leq
\frac{1}{\sqrt{G(\chi;\mu)- G(\chi-\tau(\chi-\tilde\zeta); \mu)} }\cdot \frac{\partial\chi}{\partial\omega}\leq
\frac{1}{\sqrt{- G(\tilde\zeta; \mu)} }\cdot\frac{\partial\chi}{\partial\omega}.
\end{align*}
For each $\mu\in(\mu_c, \omu)$, we can prove 
$\partial\tilde\I_1/\partial\omega$ is uniformly bounded in 
$\omega\in(\omega_*(\mu), \tilde\zeta]$, 
while in the case $\mu=\omu$ 
\[
\partial\tilde\I_1/\partial\omega \to -\infty\quad (\omega\to\omega_*(\omu)=0).
\]
This leads us to the desired assertion.
\eproof

\appendix
\section{Appendix}
We estimate $\la u(\cdot; \mu,\ell) \ra $ for $\mu\in(\mu_c, \omu]$. 
In this case $\la u(\cdot; \mu,\ell) \ra $  does not converges $0$ as $\ell\to\infty$.
As a matter of fact, we have

\begin{lem}
For the solution  
$ u(\cdot; \mu,\ell)$ obtained in \rm{(i)} and \rm{(ii)} of Lemma \ref{lem:1}, 
\begin{align*}
\lim_{\ell\to\infty}\frac{1}{\ell}\int_0^\ell u(x; \mu,\ell)\ dx
=\begin{cases}
\beta(\mu)&  (\mu \neq \omu), \\ \\
\displaystyle\frac{ \beta(\mu) }{1+\sqrt{h}}&  (\mu = \omu)
\end{cases}
\end{align*}
hold, where 
$h:=d\beta(\omu)(\beta(\omu)-\alpha(\omu))/\kappa^2(\beta(\omu)+1)$.
\end{lem}
\par\medskip
\proof
First notice
\begin{align}
\label{bd1}
\frac{1}{\ell}\int_0^\ell u(x;\mu,\ell)\ dx
\leq \frac1{\ell}\int_0^\ell \chi(\omega(\mu, \ell), \mu)dx<\beta(\mu).
\end{align}
Take any small $\ep>0$ and there is $\ell_\ep$ such that for $\ell>\ell_\ep$ 
\[
\beta(\mu)-\ep<\chi(\omega(\mu, \ell),\mu).
\]
In a similar way as in \eqref{tilrho} we have the expression 
\begin{align*}
&
\int_0^\ell u(x;\mu,\ell)\ dx
=\int_{\omega}^{\chi(\omega(\mu, \ell),\mu)}\frac{z dz}{\sqrt{2(G(\chi(\omega(\mu, \ell),\mu); \mu)-G(z;\mu))}} = \K_1+\K_2, \\
&
\K_1:=
\int_{{\beta(\mu)-\ep}}^{\chi(\omega(\mu, \ell),\mu)}\frac{z dz}{\sqrt{2(G(\chi(\omega(\mu, \ell),\mu); \mu)-G(z;\mu))}}, \\
&
\K_2:=
\int_{\omega}^{\beta(\mu)-\ep}\frac{z dz}{\sqrt{2(G(\omega(\mu,\ell); \mu)-G(z;\mu))}}.
\end{align*}
We put 
\begin{align*}
&
\ell=\K_3+\K_4,\\
&
\K_3:=
\int_{\beta(\mu)-\ep}^{\chi(\omega(\mu,\ell),\mu)}\frac{ dz}{\sqrt{2(G(\chi(\omega(\mu,\ell),\mu); \mu)-G(z;\mu))}}, \\
&
\K_4:=
\int_{\omega}^{\beta(\mu)-\ep}\frac{dz}{\sqrt{2(G( \omega(\mu,\ell); \mu)-G(z;\mu))}}.
\end{align*}
Then we have
\begin{align*}
\K_1>
(\beta(\mu)-\ep)\K_3, \qquad
\ell=\K_3+\K_4,
\end{align*}
and
\begin{align*}
\frac{1}{\ell}\int_0^\ell u(x;\mu,\ell)\ dx
=\frac{\K_1+\K_2}{\K_3+\K_4}
>\frac{(\beta(\mu)-\ep)\K_3+\K_2}{\K_3+\K_4}.
\end{align*}
We abbreviate $\omega(\mu,\ell)$ as $\omega$ below as long as no confusion.
Using the change of variable $z = \omega+ (\beta(\mu) -\ep - \omega)\sin^2\theta$, we have
\begin{align*}
	\K_2 
	& = \int_0^{\frac \pi 2} \frac{ 2\{\omega + ( \beta(\mu) -\ep -\omega )\sin^2\theta \}( \beta(\mu) -\ep -\omega )\sin\theta \cos \theta }
	{ \sqrt{2(G(\omega;\mu)-G(\omega+(\beta(\mu)-\ep-\omega)\sin^2\theta;\mu) }} d\theta
	\\
	&=  \int_0^{\frac \pi 2} \frac{ \{\omega + ( \beta(\mu) -\ep -\omega )\sin^2\theta \}\sqrt{2 (\beta(\mu) -\ep -\omega )}\cos \theta }
	{ \sqrt{ -\int_0^1 g( \omega + s( \beta(\mu) -\ep -\omega )\sin^2\theta;\mu )ds }} d\theta
	< \infty,\\
	\K_4
	& =\int_0^{\frac \pi 2} \frac{ \sqrt{2 (\beta(\mu) -\ep -\omega )} \cos\theta}
	{ \sqrt{ -\int_0^1 g( \omega + s( \beta(\mu) -\ep -\omega )\sin^2\theta ;\mu)ds }} d\theta
	< \infty,
\end{align*}
where $g(u;\mu)< 0~(0<u<\beta(\mu))$.
Next we estimate
\begin{align*}
	\K_3 &> \int_{ \beta(\mu) - \ep }^{\chi(\omega(\mu,\ell),\mu)}\frac{dz}{\sqrt{2(G(\beta(\mu);\mu) -G(z;\mu))}}  \\
	&>	
	\frac{1}{ \sqrt{C_2}} \int_{ \beta(\mu) - \ep }^{\chi(\omega(\mu,\ell),\mu) } \frac{dz}{ \beta(\mu) - z}
	= \frac{1}{ \sqrt{C_2}} \log \frac{\ep}{ \beta(\mu) - \chi(\omega(\mu,\ell),\mu)  }\\
	&\to \infty\quad (\ell \to \infty),
\end{align*}
where we put
\[
C_2 := \max_{ { {\mu_c} }\le \mu\leq\omu} \max_{\beta(\mu) -\ep \le u\le \beta(\mu)} |g_u(u; \mu)|,
\] 
and used $\lim_{\ell \to \infty} \chi(\omega(\mu,\ell),\mu) = \beta(\mu)$. 
Making use of  $\lim_{\ell\to\infty}\K_3\to\infty$ and that $\K_2, \K_4$ are bounded when $\mu\neq\omu$, we obtain for each $\mu\in(\mu_c, \omu)$
\[
\liminf_{\ell\to\infty}\frac{1}{\ell}\int_0^\ell u(x;\mu,\ell) dx
\geq \beta(\mu)-\ep.
\]
Combining this and \eqref{bd1} yields
\[
\beta(\mu)-\ep\leq \liminf_{\ell\to\infty}\frac{1}{\ell}\int_0^\ell u(x;\mu,\ell)
\leq  \limsup_{\ell\to\infty}\frac{1}{\ell}\int_0^\ell u(x;\mu,\ell)\leq \beta(\mu).
\]
Since we take $\ep$ arbitrarily, we conclude 
\[
\lim_{\ell\to\infty}\frac{1}{\ell}\int_0^\ell u(x;\mu,\ell)= \beta(\mu).
\]

As for the case 
$\mu=\omu$, we use
\begin{align*}
&
\frac{(\beta(\omu)-\ep)\K_3}{\K_3+\K_4}+\frac{\K_2}{\K_3+\K_4} \\
&\qquad
<\frac{\K_1+\K_2}{\K_3+\K_4}
<\frac{\beta(\omu)\K_3}{\K_3+\K_4}+\frac{\K_2}{\K_3+\K_4}.
\end{align*}
Since $\lim_{\ell\to\infty}\K_4=\infty$, we need to estimate
${\K_4}/{\K_3}$.  
Recall $G(u; \mu)=-u^2/2+O(|u|^3)\ (u \to 0)$.
For given $\theta\in(0,1)$ we take $\delta>0$ so that
\begin{align*}
&(1-\theta)(z^2-\omega(\omu,\ell)^2)\\
&\quad
<2(G(\omega(\omu,\ell); \omu)-G(z; \omu))<(1+\theta)(z^2-\omega(\omu,\ell)^2)\quad(0<\omega(\omu,\ell)<z<\delta).
\end{align*}
Hence, we obtain
\begin{align}
\label{K4d}
\begin{aligned}
&\int_{\omega(\omu,\ell)}^\delta \frac{dz}{\sqrt{2(G(\omega(\omu,\ell); \omu)-G(z;\omu))}}
=\int_{\omega(\omu,\ell)}^\delta \frac{dz}{\sqrt{z^2-\eta(\ell,\omu)^2}}+O(\theta)
\\
&\qquad
=
\log(\delta+\sqrt{\delta^2-\omega(\omu,\ell)^2})-
\log\omega(\omu,\ell)+O(\theta).
\end{aligned}
\end{align}

Next consider the integral near $u=\beta(\omu)$.
Since 
\[
G_u(\beta(\omu); \omu)=0,\quad
G_{uu}(\beta(\omu) ; \omu)=g_u(\beta(\omu) ; \omu)
=-\frac{d}{\kappa^2}\frac{\beta(\omu)(\beta(\omu)-\alpha(\omu))}{\beta(\omu)+1},
\]
we have
\begin{align*}
&G(u; \mu)=-\frac1{2}h(\beta(\omu)-u)^2+O(|\beta(\omu)-u|^3)\quad \ (u \to \beta(\mu)), \\
&h:=\frac{d}{\kappa^2}\cdot\frac{\beta(\omu)(\beta(\omu)-\alpha(\omu))}{\beta(\omu)+1}.
\end{align*}
Thus, for given $\theta\in(0,1)$, there is $\ep>0$ such that
\begin{align*}
&(1-\theta)h[(\beta(\omu)-z)^2-(\beta(\omu)-\oomega)^2]\\
&
<2(G(\oomega ; \omu)-G(z; \omu))<(1+\theta)h[(\beta(\omu)-z)^2-(\beta(\omu)-\oomega)^2]\\
&
 \quad(\beta(\omu)-\ep<z<\oomega<\beta(\omu)),
\end{align*}
where we simply put
\[
\ochi=\chi(\omega(\omu,\ell){,} \omu).
\]
We compute
\begin{align}
\K_3 =&\int_{\beta(\omu)-\ep}^{\ochi} \frac{dz}{\sqrt{2(G(\omega(\omu,\ell); \omu)-G(z; \omu))}} \notag\\
& =
\frac1{\sqrt{h}}\int_{\beta(\omu)-\ep}^{\ochi} \frac{dz}{\sqrt{(\beta(\omu)-z)^2-(\beta(\omu)-\ochi)^2}}+O(\theta) \notag\\
&
=\frac1{\sqrt{h}}\left[
-\log\left\{\beta(\omu)-z+\sqrt{
(\beta(\omu)-z)^2-(\beta(\omu)-\ochi)^2}\right\}
\right]_{z=\beta(\omu)-\ep}^{\ochi}+O(\theta) \notag\\
&
=\frac1{\sqrt{h}}\left[\log\left\{\ep+\sqrt{
\ep^2-(\beta(\omu)-\ochi)^2}\right\} 
- \log(\beta(\omu)-\ochi)\right]+O(\theta) \label{est:K_3}
\end{align}
In view of
\begin{align*}
&
G(\omega(\ell, \omu); \omu)=-\frac{1}{2}\omega(\ell, \omu)^2+O(\omega(\omu,\ell)^3)\ ~~~(\omega \to 0)\\
&
G(\ochi; \omu) 
=-\frac{h}{2}\{\beta(\omu)- \ochi \}^2+O(|\beta(\omu) - \ochi |^3)\ ~~~ (\chi \to \beta(\mu))
\end{align*}
and $G(\omega(\ell, \omu); \omu)=G( \ochi ; \omu)$, we see
\begin{align}
\label{xieta}
\sqrt{\chi}\{\beta(\omu)- \ochi \}=\omega(\omu,\ell)+O(\omega(\omu,\ell)^{3/2}).
\end{align}

Taking $\delta=\ep$ in \eqref{K4d}, we have
\begin{align}
\label{K4a}
\begin{aligned}
\K_4&=\log(\ep+\sqrt{\ep^2-\omega(\omu,\ell)^2})-
\log\omega(\omu,\ell)+O(\theta)\\
& \quad
+\int_{\ep}^{\beta(\omu)-\ep} \frac{dz}{\sqrt{2(G(\omega(\omu,\ell); \omu)-G(z;\omu))}}.
\end{aligned}
\end{align}
By \eqref{xieta} we get 
\begin{align*}
&
\log(\ep+\sqrt{\ep^2-\omega(\omu,\ell)^2})-
\log\omega(\omu,\ell) \\
&
=\log[\ep+\sqrt{\ep^2-\chi\{\beta(\omu)- \ochi \}^2}]-
\log[\sqrt\chi\{\beta(\omu) - \ochi \}],
\end{align*}
where, if necessary, we take $\ell$ larger so that 
\[
\sqrt\chi\{\beta(\omu)- \ochi \}<\ep
\]
holds.
Applying 
\begin{align*}
&\frac{\log \left[ \ep+\sqrt{\ep^2-\chi\{\beta(\omu) - \ochi \}^2} \right]-
\log ( \sqrt\chi\{\beta(\omu) - \ochi \}) }
{\frac1{\sqrt{h}}\left[\log\left\{\ep+\sqrt{
\ep^2-\{\beta(\omu) - \ochi \}^2}\right\} 
- \log(\beta(\omu) - \ochi )\right]}\to\sqrt{h} \\
& (\ell\to\infty),
\end{align*}
by $\omega(\ell, \omu) \to \omega_*(\omu)=0, \ \ochi \to \beta(\omu)$ as $\ell \to \infty$ due to Lemma \ref{lem:trhosol}
and that the integral in the right-hand side of \eqref{K4a}
is bounded as $\ell\to\infty$, we obtain
\begin{align*}
\lim_{\ell\to\infty}\frac{\K_4}{\K_3}=\lim_{\ell\to\infty}
\frac{\sqrt{h}\log\omega(\omu,\ell)}{\log\{\beta(\omu) - \ochi  \}}
=\sqrt{h}
\end{align*}
where the first equality is given by using \eqref{est:K_3}, \eqref{K4a} and L'H\^{o}pital's rule.

\par\bigskip\par


Consequently, we obtain
\begin{align*}
&
\lim_{\ell\to\infty}\left(\frac{(\beta(\omu)-\ep)\K_3}{\K_3+\K_4}+\frac{\K_2}{\K_3+\K_4}\right)=\frac{\beta(\omu)-\ep}{1+
\sqrt{h}},
 \\
 &
 \lim_{\ell\to\infty}\left(\frac{\beta(\omu)\K_3}{\K_3+\K_4}
 +\frac{\K_2}{\K_3+\K_4}\right)=\frac{\beta(\omu)}{1+\sqrt{h}}.
 \end{align*}
 Since $\ep$ is arbitrarily taken, we have the assertion
 \[
 \lim_{\ell\to\infty} \frac{1}{\ell}\int_0^\ell \tilde{u}(x; \omu, \ell)\ dx
 =\frac{\beta(\omu)}{1+\sqrt{h}},
 \]
 for the case $\mu=\omu$.
\eproof

%
%
%

\par\vspace{0.5cm}

\section*{Acknowledgments}
The authors would like to express their sincere acknowledgement to Professor Yasumasa Nishiura for letting know them the article of Yochelis-Beta-Giv. 
The authors were partially supported by JSPS KAKENHI Grant Number  22K03444 and the second author was by JSPS KAKENHI Grant Numbers 20K14364.

\end{document}